\documentclass[10pt,leqno]{article}
\usepackage{amsmath,amsfonts,amscd,amsthm}
\usepackage{graphicx}
\usepackage{multicol}
\usepackage{slashbox}
\usepackage{cite}
\newcommand{\R}{\mathbb{R}}
\newcommand{\N}{\mathbb{N}}
\newcommand{\dps}{\displaystyle}
\newcommand{\bu}{{\bf u}}
\newcommand{\uhb}{\bu_{h,{\beta}}}

\newcommand{\td}{{\tilde\Delta}}
\newcommand{\pr}{\partial}
\newcommand{\pt}{\partial_t}
\newcommand{\bw}{{\bf w}}
\newcommand{\e}{{\bf e}}

\newcommand{\bJ}{{\bf J}}

\newcommand{\bU}{{\bf U}}
\newcommand{\tbU}{\tilde{\bU}}

\newcommand{\bv}{{\bf v}}

\newcommand{\bta}{\mbox{\boldmath $\eta$}}
\newcommand{\hbta}{\widehat{\bta}}
\newcommand{\tbta}{\tilde{\bta}}

\newcommand{\bphi}{\mbox{\boldmath $\phi$}}

\newcommand{\bxi}{\mbox{\boldmath $\xi$}}
\newcommand{\tbxi}{\tilde{\bxi}}
\newcommand{\hbxi}{\widehat{\bxi}}
\newcommand{\bV}{{\bf V}}
\newcommand{\bW}{{\bf W}}
\newcommand{\bH}{{\bf H}}
\newcommand{\bL}{{\bf L}}
\newcommand{\f}{{\bf f}}
\newcommand{\tf}{\tilde{\f}}
\newcommand{\g}{{\bf g}}
\newcommand{\x}{{\bf x}}

\newcommand{\ve}{\varepsilon}

\begin{document}

\newcommand{\se}{\setcounter{equation}{0}}
\def\theequation{\thesection.\arabic{equation}}

\newtheorem{theorem}{Theorem}[section]
\newtheorem{cdf}{Corollary}[section]
\newtheorem{lemma}{Lemma}[section]
\newtheorem{remark}{Remark}[section]
\newtheorem{example}{Example}[section]

\def\cydot{\leavevmode\raise.4ex\hbox{.}}

\title
{Backward Euler method for the equations of motion arising in Oldroyd model of order one with nonsmooth initial data}
\author{Bikram Bir\thanks{Department of Mathematical Sciences,
Tezpur University, Napaam, Sonitpur, Assam-784028, India. Email: bikramb@tezu.ernet.in, deepjyoti@tezu.ernet.in}, ~ Deepjyoti Goswami\footnotemark[1] ~~and
Amiya K. Pani\thanks{Department of Mathematics, 
Indian Institute of Technology Bombay, Powai, Mumbai-400076, India. Email: akp@math.iitb.ac.in}
}
\date{}
\maketitle

\begin{abstract}
In this paper, a backward Euler method combined with finite element discretization in spatial direction 
is discussed for the equations of motion arising in 
the $2D$ Oldroyd model of viscoelastic fluids of order one with the forcing term independent 
of time or in $L^{\infty}$ in time. It is shown that the estimates of the discrete solution 
in Dirichlet norm is bounded uniformly in time. Optimal {\it a priori} error estimate in 
$\bL^2$-norm is derived for the discrete problem with non-smooth initial data. This estimate 
is shown to be uniform in time, under the assumption of uniqueness condition. Finally, we
present some numerical results to validate our theoretical results.
\end{abstract}

\vspace{1em} 

\section{Introduction}
In this paper, we consider fully-discrete approximations to the equations of motion arising in the Oldroyd fluids (see Oldroyd \cite{Old}) of order one: 
\begin{eqnarray}\label{om}
~~\frac {\partial \bu}{\partial t}+\bu\cdot\nabla\bu-\mu\Delta\bu-\int_0^t \beta
 (t-\tau)\Delta\bu (\tau)\,d\tau+\nabla p=\f, ~~~~~\mbox {in}~ \Omega,~t>0
\end{eqnarray}
with incompressibility condition
\begin{eqnarray}\label{ic}
 \nabla \cdot \bu=0,~~~~~\mbox {on}~\Omega,~t>0,
\end{eqnarray}
and initial and boundary conditions
\begin{eqnarray}\label{ibc}
 \bu(x,0)=\bu_0~~\mbox {in}~\Omega,~~~\bu=0,~~\mbox {on}~\partial\Omega,~t\ge 0.
\end{eqnarray}
Here, $\Omega$ is a bounded domain in $\mathbb{R}^2$ with boundary $\partial \Omega$,
$\mu = 2 \kappa\lambda^{-1}>0$, the kernel $\beta (t) = \gamma \exp (-\delta t),~
\gamma= 2\lambda^{-1}(\nu-\kappa \lambda^{- 1})>0$ and $\delta =\lambda^{-1}>0$, where $\nu>0$ is
the kinematic coefficient of viscosity, $\lambda>0$ is the relaxation time and $\kappa>0$ is the retardation time. Unknowns $\bu$ and $p$ represent the velocity and the pressure of the fluid, respectively. Further, the forcing term $\f$ and initial velocity $\bu_0$ are given functions in their respective domains of definition. For more details on the model, we refer to \cite{Old}. 

The model has been studied for more than three decades now; for early works and, for a brief introduction on the continuous and semi-discrete cases, we refer to \cite{HLSST,GP11,PY05} and references, therein. For recent and other notable works, we refer to \cite{GH16,GH18,ZQ18,ZQJY18,ZY15,ZZQ18, AD19, LS19, MO20, MO20_1, YLS20} and references, therein.

Our present investigation is a continuation of our works in \cite{GP11} and \cite{PYD06}. In \cite{GP11}, {\it a priori} estimates and regularity results for the solution pair $\{\bu,p\}$ of (\ref{om})-(\ref{ibc}) are established under realistically assumed data and when $\f,\f_t\in L^{\infty}(\R_{+};\bL^2(\Omega))$. Further, optimal error estimates for the velocity and for the pressure are derived for the semidiscrete Galerkin approximations and these results are shown to be uniform under uniqueness condition. In this work, a completely discrete scheme based on backward Euler method is developed and analyzed for the problem (\ref{om})-(\ref{ibc}). 
With $k$ as uniform time-step size, $t_n=nk$ as $n$th time level and $t_N=Nk$ as the final time,  we denote $\bU^n$ as the approximation of semi-discrete solution $\bu_h$ at $t=t_n$ and fully-discrete approximation of $\bu$ at $t=t_n$. We analyze the error due to the approximation $\bU^n$ in light of non-smooth initial data
and present optimal order error estimates. Before we discuss our main result and highlight the technical difficulties in its proof, we would first like to have a look at the available literature in this direction.

Literature for the fully-discrete approximation to the problem (\ref{om})-(\ref{ibc}) is limited. In \cite{AO}, Akhmatov and Oskolkov have discussed stable and convergent finite difference schemes for the problem (\ref{om})-(\ref{ibc}). On the other hand, Pani {\it et al.}, in \cite{PYD06}, have considered a linearized backward Euler method to discretize in time direction only keeping spatial direction continuous and have used semi-group theoretic approach to establish {\it a priori} error estimates.  The following time discrete error bounds are proved in \cite{PYD06}:
(for $0<\alpha<\min\{\delta,\lambda_1\}$, $\lambda_1$ is the smallest positive eigenvalue of the Stokes operator)
$$ \|\bu(t_n)-\bU^n\| \le Ce^{-\alpha t_n}k,~~ 
\|\bu(t_n)-\bU^n\|_1 \le Ce^{-\alpha t_n}k(t_n^{-1/2}+\log\frac{1}{k}), $$
for smooth initial data, i.e., $\bu_0\in\bH^2(\Omega)\cap\bJ_1$ (see, $\bJ_1$ in Section 2) and for zero forcing term ($\f\equiv 0$).
 In \cite{WHS10}, Wang {\it et al.} have extended this work to non-zero forcing function.
They have used energy arguments, along with uniqueness condition 
to obtain for fully discrete solution $\bU^n$, the following uniform error estimates 
$$ \|\bu(t_n)-\bU^n\| \le C(h^2+k),~~
(\tau^*)^{1/2}\|\bu(t_n)-\bU^n\|_1 \le C(h+k), $$
where $\tau^*(t_n)= \min\{1,t_n\}$ again for smooth data.
In \cite{GH16}, Guo {\it et al.} have worked with a second-order time discretization scheme based on Crank-Nicolson/Adams-Bashforth as part of fully discrete analysis and have derived optimal error estimate  under smooth initial data.

In the present work, we examine on the backward Euler method and 
prove the following (see, Theorem \ref{l2eebe}), when $\bu_0\in\bH_0^1(\Omega)$ :
\begin{align*}
\|\bu_h(t_n)-\bU^n\| \le K_nt_n^{-1/2}k\big(1+\log\frac{1}{k}\big)^{1/2},~~1\le n\le N <  +\infty,
\end{align*}
where $K_n>0$ is the error constant that depends only on the given data and, in particular, is
independent of both $h$ and $k$. But it grows exponentially with time, that is, $K_n \sim
O(e^{t_n})$ and therefore, the above error estimate is local (in time). Under uniqueness
condition, we have shown the error to be uniformly bounded as $t\to +\infty$, see Section $6$.

These results are proved for nonsmooth initial data, that is, under realistically assumed regularity on the exact solution of the problem (\ref{om})-(\ref{ibc}).
For example, Lemma \ref{dth2} says that $\|\bu_h(t)\|_2$ and $\|\bu_{ht}\|$ are of $O(t^{-1/2})$.  
As in \cite{GP11}, this breakdown at $t=0$ is a major bottle-neck in our error analysis. To illustrate our point, consider $\bu_0\in \bH_0^1(\Omega)\cap\bH^2(\Omega)$ (smooth initial data). Then, the error $\e_n= \bU^n-\bu_h(t_n)$ satisfies the following estimate (see, Lemma $4.2$ of \cite{WHS10}):
$$ \|\e_n\| \sim O(k),~~1\le n \le N, $$
Following similar argument but now with $\bu_0\in \bH_0^1(\Omega)$, we would only obtain
$$ \|\e_n\| \sim O(k^{1/2}\big(1+\log\frac{1}{k}\big)^{1/2}),~~1\le n \le N. $$
The loss in the order of $k$, in fact, is due to the singularity of the higher-order norms of the solution at $t=0$.
The standard technique, in such cases, is to multiply by a weight $t^r,~r\in \N$ to compensate for this singularity, thereby, recovering full order of convergence. But in our case, a direct application of this technique fails due to the presence of the memory term. Note that the kernel $\beta$ present in the equation (\ref{om}) has a certain positivity property (see Lemma $2.1$ of \cite{GP11}) and we choose our quadrature rule to conform with it, see (\ref{rrp}).
This is crucial to our analysis. 
But when we opt for weighted Sobolev norm with a weight $t^r,~r\in \N$, it nullifies the positivity property of our quadrature rule.
The main effort of this work is to overcome this difficulty and to recover optimal
fully discrete error estimate for the velocity.
%
This requires borrowing certain tools from error analysis of linear parabolic integro-differential equations  with non-smooth data  (see; \cite{PS98, PS198, TZ89}), like, the summation technique (we call it here "hat operator", see (\ref{sum0})), which adds to the technicality.

The error analysis has been carried out by splitting the error into two parts such as error due to linearized part and error due to nonlinear part and then analyzing both part separately. 
Our approach has two advantages. Firstly, we notice that the exponential increase of the error bound (as $t\to +\infty$) is due to the nonlinear part, since the other part is uniformly bounded as $t\to +\infty$, see, lemma \ref{l2eebxi}. Secondly, for non-smooth initial data, the uniform error estimate (which has been done in Section $6$), can only be obtained by splitting the error likewise.

Apart from the error analysis, we have also established uniform (in time) Dirichlet norm
for the fully discrete solution $\bU^n$, meaning $\|\nabla\bU^n\|,~1\le n\le N$ remains 
bounded as $t_N\to +\infty$. It is crucial for the long-time stability of the implicit scheme. In 
case of Navier-Stokes, the proof of the Dirichlet norm of $\bU^n$, which is valid for all 
time, involves applying discrete version of the uniform Gronwall's Lemma (see, Lemma 2.6 from 
\cite{TW06}). In our case, it is difficult to apply the uniform Gronwall's Lemma due to the 
presence of the quadrature term. 
Hence, we resort to a new way of looking into the ideas behind the uniform Gronwall's Lemma to 
establish our result.

\noindent We now summarize our main results as follows:
\begin{itemize}
\item [(i)] Uniform bound in time for the fully discrete solution in the Dirichlet norm depicting long term stability.
\item [(ii)] New uniform estimates for the error
associated with fully discrete linearized problem.
\item [(iii)] Local optimal error estimates for the discrete velocity in $\bL^2$-norm.
\item [(iv)] Optimal global fully discrete error estimates
under the uniqueness assumption.

\end{itemize}

At this stage, it is useful to compare our results with the results derived in \cite{PYD06}. In \cite{PYD06}, only discretization in time keeping spatial variables continuous has been analyzed using semigroup theory approach for the homogeneous problem, that is, $\f=0$, and error estimates are obtained under the assumption of $\bu_0\in\bH^2\cap\bJ_1$. But the present analysis deals with the fully discrete scheme, uses energy arguments and establishes the optimal error estimates for the nonhomogeneous problem (\ref{om})-(\ref{ibc}) with nonsmooth initial data, that is, $\bu_0\in\bH_0^1\cap\bJ_1$. In both papers, a common thread is the time weighted estimates.

\noindent
Henceforward, we will use $K$ and $C$ as positive generic constants, where $K$ would depend on the given data.

The remaining part of this paper is organized as follows. In Section $2$, we discuss some notations, basic assumptions and  weak formulations of the problem (\ref{om})-(\ref{ibc}).  Section $3$ deals with a brief description of a semi-discrete Galerkin finite element method. Section $4$ is devoted to backward Euler method. Optimal error bounds are obtained for the velocity and for the pressure for the problem (\ref{om})-(\ref{ibc}) with non-smooth initial data, in Section $5$, whereas, in Section $6$, these bounds are shown to be uniform in time, under uniqueness condition. In Section $7$, we present some numerical experiments whose results confirm our theoretical findings. Finally, we summarize our results in  Section $8$.

\section{Preliminaries}
\se

For our subsequent use we denote by bold face letters the $\R^2$-valued
function space such as
\begin{align*}
\bH_0^1 = [H_0^1(\Omega)]^2,~~~\bL^2 = [L^2(\Omega)]^2
 ~~\mbox{and}~~~\bH^m=[H^m(\Omega)]^2,
\end{align*}
where $H^m(\Omega)$ is the standard Hilbert Sobolev space of order $m$. Note
that $\bH^1_0$ is equipped with a norm
$$ \|\nabla\bv\|= \left(\sum_{i,j=1}^{2}(\partial_j v_i, \partial_j
 v_i)\right)^{1/2}=\left(\sum_{i=1}^{2}(\nabla v_i, \nabla v_i)\right)^{1/2}. $$
Further, we introduce below, divergence free function spaces:
\begin{align*}
\bJ_1 &= \{\bphi\in\bH_0^1 : \nabla \cdot \bphi = 0\}, \\
\bJ &= \{\bphi \in\bL^2 :\nabla \cdot \bphi = 0~\mbox{in}~
 \Omega,~\bphi\cdot{\bf n}|_{\pr \Omega}=0~~\mbox{holds weakly}\}, 
\end{align*}
where ${\bf n}$ is the outward normal to the boundary $\pr\Omega$ and $\bphi
\cdot {\bf n} |_{\pr\Omega} = 0$ should be understood in the sense of trace
in $\bH^{-1/2}(\pr\Omega)$, see \cite{temam}. Let $H^m/\R$ be the
quotient space consisting of equivalence classes of elements of $H^m$ differing
by constants, which is equipped with norm $\| p\|_{H^m /\R}=\inf_{c\in\R} \| p+c\|_m$.
For any Banach space $X$, let $L^p(0, T; X)$ denote the space of measurable $X$-valued functions $\bphi$ on $(0,T),~0<T\le +\infty$, such that
$$ \int_0^T \|\bphi (t)\|^p_X~dt <+\infty~~~\mbox{if}~~1 \le p < +\infty, $$
and for $p=+\infty$
$$ {\dps{ess \sup_{0<t<T}}} \|\bphi (t)\|_X <+\infty. $$
Through out this paper, we make the following assumptions: \\
(${\bf A1}$). For $\g\in \bL^2(\Omega)$, let the unique pair of solutions $\{\bv,q\}\in 
\bJ_1\times L^2(\Omega)/\R $ for the steady state Stokes problem
\begin{align*}
-\Delta {\bv} + \nabla q = {\bf {g}},~
\nabla \cdot\bv = 0~{\mbox in}~\Omega~~\mbox{and}~~\bv|_{\pr \Omega}=0,
\end{align*}
satisfy the following regularity result
$$  \|\bv\|_2 + \|q\|_{H^1(\Omega)/\R} \le C\|\g\|. $$
\noindent
(${\bf A2}$). The initial velocity $\bu_0$ and the external force $\f$ satisfy, for
positive constant $M_0,$
$$ \bu_0\in\bJ_1,~\f,\f_t \in L^{\infty}(\R_{+};\bL^2(\Omega))~~\mbox{with}~~ \|\bu_0\|_1 \le M_0,~~{\dps{\sup_{t>0} }}\big\{\|\f\|, \|\f_t\|\big\} \le M_0. $$

\noindent With $P$ as orthogonal projection of $\bL^2(\Omega)$ onto $\bJ$, we set $\td= P\Delta$ as the Stokes operator. Then, (${\bf A1}$) implies
\begin{align*}
\|\bv\|_2 \le C\|\td\bv\|,~~\|\nabla\bv\|^2 \le \lambda_1^{-1} \|\td\bv\|^2,~~\forall ~\bv\in\bJ_1\cap\bH^2(\Omega), \mbox{  and ~}
\|\bv\|^2 \le \lambda_1^{-1} \|\nabla\bv\|^2,~~\forall ~\bv\in\bJ_1,
\end{align*}
where  $\lambda_1$ is the least positive eigenvalue of the Stokes operator, see \cite{HR82}.

Before going to the details, let us introduce the weak formulation of (\ref{om})-(\ref{ibc}): Find a pair of functions $\{\bu(t), p(t)\}\in \bH_0^1\times L^2/\mathbb{R},~t>0,$ such that
\begin{align*}
(\bu_{t}, \bphi) +\mu  (\nabla\bu,\nabla\bphi)+   (\bu\cdot\nabla\bu,\bphi)+ \int_0^t\beta(t-s)(\nabla\bu(s),\nabla\bphi)ds -(p, \nabla \cdot \bphi) &=(\f, \bphi),~~\forall \bphi\in\bH_0^1 \\
(\nabla \cdot \bu, \chi) &=0,~~\forall\chi\in L^2.	 \nonumber
\end{align*}
Equivalently, find $\bu(\cdot,t)\in \bJ_1$ such that 
\begin{align*}
(\bu_{t}, \bphi) +\mu  (\nabla\bu,\nabla\bphi)+   (\bu\cdot\nabla\bu,\bphi)+ \int_0^t\beta(t-s)(\nabla\bu(s),\nabla\bphi)ds =(\f, \bphi),~~\forall \bphi\in\bJ_1, ~ t>0.
\end{align*}

\noindent Now, we present below, the discrete Gronwall's Lemma. For a proof, see, \cite{HL98,PTW92}: %
\begin{lemma}[discrete Gronwall's Lemma]\label{gl}
Let $k,B$ and $\{a_i,b_i,c_i,d_i\}_{i\in\mathbb{N}}$ be non-negative numbers such that
\begin{equation*}
a_n+k\sum_{i=1}^n b_i \le B+k\sum_{i=1}^{n-1} c_i+k\sum_{i=1}^{n-1} d_ia_i, ~~\forall n\ge 1.
\end{equation*}
Then,
\begin{equation*}
a_n+k\sum_{i=1}^n b_i \le \Big\{B+k\sum_{i=1}^{n-1} c_i\Big\}\exp\Big(k\sum_{i=1}^{n-1} 
d_i\Big), ~~\forall n\ge 1.
\end{equation*}
\end{lemma}
%
%

\section{Semidiscrete  Galerkin Approximations}
\se
From now on, we denote $h$ with $0<h<1$ to be a real positive discretization
parameter tending to zero. Let  $\bH_h$ and $L_h$ be
finite dimensional subspaces of $\bH_0^1 $ and $L^2$, respectively,
approximating the velocity vector and the pressure. Assume the following
approximation properties for the spaces $\bH_h$ and $L_h$: \\
${\bf (B1)}$ For each $\bw \in \bH_0^1 \cap \bH^2 $ and $ q \in
H^1/\R$ there exist approximations $i_h w \in \bH_h $ and $ j_h q \in
L_h $ such that
\[ \|\bw-i_h\bw\|+ h \| \nabla (\bw-i_h \bw)\| \le Ch^2 \| \bw\|_2,
 ~~~~\| q - j_h q \|_{L^2/\R} \le Ch \| q\|_{H^1/\R}. \]
Further, suppose that the following inverse hypothesis holds for $\bw_h\in\bH_h$:
\begin{align*}
 \|\nabla \bw_h\| \leq  Ch^{-1} \|\bw_h\|.
\end{align*}
To define the Galerkin approximations, set for $\bv, \bw, \bphi \in \bH_0^1$,
$$ a(\bv, \bphi) = (\nabla \bv, \nabla \bphi) ~~~\mbox{and}~~~
b(\bv,\bw,\bphi)= \frac{1}{2} (\bv \cdot \nabla \bw , \bphi)
   - \frac{1}{2} (\bv \cdot \nabla \bphi, \bw). $$
Note that the operator $b(\cdot, \cdot, \cdot)$ preserves the antisymmetric property of
the original nonlinear term, that is,
$$ b(\bv_h, \bw_h, \bw_h) = 0, \;\;\; \forall \bv_h, \bw_h \in {\bH}_h. $$
Now the semidiscrete Galerkin formulation reads as:  Find $(\bu_h(t), p_h(t))\in
\bH_h\times L_h$ such that $ \bu_h(0)= \bu_{0h} $ 
\begin{align}\label{dwfh}
(\bu_{ht}, \bphi_h) +\mu a (\bu_h,\bphi_h)+& b(\bu_h,\bu_h,\bphi_h)+ a(\uhb, 
\bphi_h) -(p_h, \nabla \cdot \bphi_h) =(\f, \bphi_h), \\
&(\nabla \cdot \bu_h, \chi_h) =0, \nonumber
\end{align}
for $\bphi_h\in\bH_h,~\chi_h \in L_h$ and for $t>0$. Here $\bu_{0h} \in\bH_h $ is a suitable 
approximation of $\bu_0\in \bJ_1$ and
\begin{equation*}
\uhb(t)=\int_0^t \beta(t-s) \bu_h(s)~ds.
\end{equation*}

\noindent In order to consider a discrete space analogous to $\bJ_1$, we
impose the discrete incompressibility condition on $\bH_h$ and call it as
$\bJ_h$. We define $\bJ_h$ as
$$ \bJ_h := \{ v_h \in\bH_h : (\chi_h,\nabla\cdot v_h)=0,~~\forall \chi_h \in L_h \}. $$
Note that $\bJ_h$ is not a subspace of $\bJ_1$. With $\bJ_h$ as above, we now introduce
an equivalent Galerkin formulation. Find $\bu_h(t)\in \bJ_h $ such that $\bu_h(0) =
\bu_{0h} $ and for $t>0$
\begin{equation}\label{dwfj}
(\bu_{ht},\bphi_h) +\mu a (\bu_h,\bphi_h) + a(\uhb, \bphi_h)
= -b( \bu_h, \bu_h, \bphi_h)+(\f,\bphi_h),~~\forall \bphi_h \in \bJ_h.
\end{equation}
Since $\bJ_h$ is finite dimensional, the problem (\ref{dwfj}) leads to a system of
nonlinear integro-differential equations. For global existence of a unique solution $\bu_h$
of (\ref{dwfj}), we refer to \cite{PY05}. Then $\ p_h$ is recovered from (\ref{dwfh}). Uniqueness (of $p_h$)
is obtained in the quotient space $L_h/N_h$, where
$$ N_h=\{q_h\in L_h:(q_h, \nabla\cdot\bphi_h)=0,~~\forall \bphi_h\in\bH_h\}. $$
The norm on $ L_h/N_h $ is given by
$$
 \| q_h\|_{L^2/N_h} = {\dps{\inf_{\chi_h \in N_h} }} \| q_h + \chi_h\|. $$
For continuous dependence of the discrete pressure $p_h (t) \in L_h/N_h$ on the
discrete velocity $\bu_h(t) \in \bJ_h$, we assume the following discrete
inf-sup (LBB) condition for the finite dimensional spaces $\bH_h$ and $L_h$:\\
\noindent
${\bf (B2')}$  For every $q_h \in L_h$, there exists a non-trivial function
$\bphi_h \in\bH_h$ such that
$$ |(q_h, \nabla\cdot \bphi_h)| \ge C\|\nabla \bphi_h \|\| q_h\|_{L^2/N_h}. $$
Moreover, we also assume that the following approximation property holds true
for $\bJ_h $. \\
\noindent
${\bf (B2)}$ For every $\bw \in \bJ_1 \cap\bH^2, $ there exists an
approximation $r_h \bw \in {\bf J_h}$ such that
$$ \|\bw-r_h\bw\|+h \| \nabla (\bw - r_h \bw) \| \le Ch^2 \|\bw\|_2 . $$
This is a less restrictive condition than (${\bf B2'}$) and it has been used to
derive the following properties of the $L^2$ projection $P_h:\bL^2\mapsto \bJ_h$.
We state below these results without proof. For a proof, see \cite{HR82}.
\begin{align*}
 \|\bphi- P_h \bphi\|+ h \|\nabla P_h \bphi\| &\leq Ch\|\nabla \bphi\|,~~\text{for}~ \bphi\in \bJ_h,  \\
 \|\bphi-P_h\bphi\|+h\|\nabla(\bphi-P_h \bphi)\| &\le Ch^2\|\td\bphi\|,~~\text{for}~\bphi \in \bJ_1 \cap \bH^2.
 %
\end{align*}
We now define the discrete Laplace operator $\Delta_h: \bH_h \mapsto \bH_h$ through the
bilinear form $a (\cdot, \cdot)$ as
\begin{eqnarray*}
 a(\bv_h, \bphi_h) = (-\Delta_h\bv_h, \bphi),~~\forall \bv_h, \bphi_h\in\bH_h.
\end{eqnarray*}
Denoting the orthogonal projection of $L^2$ onto $J$ as $P$, we define the Stokes operator as $\td
=P \Delta $ and set its discrete version as $\td_h = P_h \Delta_h $. The restriction of $\td_h$ to
$\bJ_h$ is invertible and we denote the inverse by $\td_h^{-1}$. Since $-\td_h$ is self-adjoint and 
positive definite, we define {\it discrete Sobolev norms} on $\bJ_h$ as follows:
$$ \|\bv_h\|_r = \|(-\td_h)^{r/2}\bv_h\|,~~\bv_h\in\bJ_h,~r\in\R. $$
We note that in particular $\|\bv_h\|_0=\|\bv_h\|$ and $\|\bv_h\|_1=\|\nabla\bv_h\|$ for $\bv_h\in\bJ_h$, and $\|\cdot\|_2$ and $\|\td_h\cdot\|$ are equivalent norms on $\bJ_h$. For further detail, we refer to \cite{HR82, HR90}.
\begin{remark}
To avoid confusion as to whether $\|\cdot\|_1$ means standard or discrete Sobolev norm, we follow the convention that if $\bv$ belongs to $\bJ_h$ then $\|\bv\|_1$ represents $\bv$ in discrete Soboelev norm, otherwise it is the standard Sobolev norm.
\end{remark}

\noindent
Below we present some estimates of the nonlinear operator $b$ for our subsequent use. The proofs of these estimates are well known and can be found in the literature based on Navier-Stokes equations (e.g., see \cite[(3.7)]{HR90}).
\begin{lemma}\label{nonlin}
Suppose conditions (${\bf A1}$), (${\bf B1}$) and (${\bf B2}$) are satisfied. Then there 
exists a positive constant $C$ such that for $\bv,\bw,\bphi\in\bH_h$, the following holds:
\begin{equation}\label{nonlin1}
 |(\bv\cdot\nabla\bw,\bphi)| \le C \left\{
\begin{array}{l}
 \|\bv\|^{1/2}\|\nabla\bv\|^{1/2}\|\nabla\bw\|^{1/2}\|\Delta_h\bw\|^{1/2}
 \|\bphi\|, \\
 \|\bv\|^{1/2}\|\Delta_h\bv\|^{1/2}\|\nabla\bw\|\|\bphi\|, \\
 \|\bv\|^{1/2}\|\nabla\bv\|^{1/2}\|\nabla\bw\|\|\bphi\|^{1/2}
 \|\nabla\bphi\|^{1/2}, \\
 \|\bv\|\|\nabla\bw\|\|\bphi\|^{1/2}\|\Delta_h\bphi\|^{1/2}, \\
 \|\bv\|\|\nabla\bw\|^{1/2}\|\Delta_h\bw\|^{1/2}\|\bphi\|^{1/2}
 \|\nabla\bphi\|^{1/2}.
\end{array}\right.
\end{equation}
\end{lemma}

\noindent Examples of subspaces $\bH_h$ and $L_h$ satisfying assumptions (${\bf B1}$),
(${\bf B2}'$), and (${\bf B2}$) can be found in \cite{BF, BP, GR}. \\
We present below a Lemma that deals with higher order estimates of $\bu_h,$ which will be
useful in the error analysis of backward Euler method for non-smooth data.
\begin{lemma}\label{dth2}
Let $0<\alpha <\min\{\delta,\lambda_1\mu\}$ and let (${\bf A1}$), (${\bf A2}$), (${\bf B1}$) and
(${\bf B2}$) be satisfied. Moreover, let $\bu_h(0)\in\bJ_h$. Then $\bu_h,$ the solutions of the
semidiscrete Oldroyd problem (\ref{dwfj}), satisfies the following {\it a priori} estimates:
\begin{align}
 \tau^*\|\bu_h\|_2^2+(\tau^*)^{r+1}(t)\|\bu_{ht}\|_r^2 & \le K,~~~~r\in \{0,1\}, 
\label{dth11} \\
 e^{-2\alpha t}\int_0^t e^{2\alpha s}(\tau^*)^r(s)\|\bu_{hs}\|^2_r\,ds & \le K, ~~~~r\in 
\{0,1,2\}, \label{dth12} \\
  e^{-2\alpha t} \int_0^t e^{2\alpha s}(\tau^*)^{r+1}(s)\|\bu_{hss}\|_{r-1}^2~ds & \le 
 K,~~~~r\in \{-1,0,1\},  \label{dth13}
\end{align}
where $\tau^*(t)=\min\{1,t\}$ and $K$ depends on the given data, but not on time $t$.
\end{lemma}
 
\begin{proof}
 The estimates (\ref{dth11})-(\ref{dth12}) can be proved as in the continuous case, see 
\cite{GP11}. The final estimate follows in a similar manner. For the sake of completeness,
we provide a sketch here. Differentiate (\ref{dwfj}) to find that, for $\bphi_h \in \bJ_h,$
\begin{align}\label{dwfjt}
(\bu_{htt},\bphi_h) +\mu a(\bu_{ht},\bphi_h) & + \beta(0)a(\bu_h,\bphi_h)-\delta
 \int_0^t \beta(t-s) a(\bu_h(s), \bphi_h)~ds \nonumber \\
 & = -b(\bu_{ht},\bu_h,\bphi_h)-b(\bu_h,\bu_{ht},\bphi_h)+(\f_t,\bphi_h).
\end{align}
Put $\bphi_h=(\tau^*)^2(t)e^{2\alpha t}\bu_{htt}$ in (\ref{dwfjt}) to obtain
\begin{align}\label{dwfjt1}
(\tau^*)^2(t)e^{2\alpha t}\|\bu_{htt}\|^2 +\mu a(\bu_{ht},(\tau^*)^2(t)e^{2\alpha t}\bu_{htt}) = (\tau^*)^2(t)e^{2\alpha t}\Big( -\gamma a(\bu_h,\bu_{htt})+\delta
 \int_0^t \beta(t-s) a(\bu_h(s),  \bu_{htt})~ds \nonumber \\
   -b(\bu_{ht},\bu_h, \bu_{htt})-b(\bu_h,\bu_{ht},\bu_{htt})+(\f_t,\bu_{htt})\Big).
\end{align}
 The second term on the left hand side of (\ref{dwfjt1}) can be written as
\begin{align*}
\mu\ a(\bu_{ht},(\tau^*)^2(t)e^{2\alpha t}\bu_{htt}) &=\frac{\mu}{2}(\tau^*)^2(t)e^{2\alpha t}
\frac{d}{dt}\|\bu_{ht}\|_1^2 \\
&= \frac{\mu}{2}\frac{d}{dt}\big( (\tau^*)^2(t) e^{2 \alpha t}
\|\bu_{ht}\|_1^2\big)-\mu\big(\alpha  (\tau^*)^2(t)+\tau^*(t)\frac{d}{dt}(\tau^*)(t)\big) e^{2\alpha t} \|\bu_{ht}\|_1^2.
\end{align*}
 Next use the Cauchy-Schwarz inequality in first, second and last terms on the right hand side of (\ref{dwfjt1}) to obtain
\begin{align}\label{dth001}
(\tau^*)^2(t) e^{2\alpha t}\|\bu_{htt}\|^2+&\frac{\mu}{2}\frac{d}{dt}\big( (\tau^*)^2(t) 
e^{2 \alpha t}\|\bu_{ht}\|_1^2\big) \le \big(\alpha (\tau^*)^2(t)+\tau^*(t)\big) 
e^{2\alpha t} \|\bu_{ht}\|_1^2 \nonumber \\
& +\gamma (\tau^*)^2(t)e^{2\alpha t}\|\bu_h\|_2\|\bu_{htt}\|+\delta (\tau^*)^2(t) 
e^{2\alpha t}\int_0^t \beta(t-s) \|\bu_h(s)\|_2\|\bu_{htt}\|~ds \nonumber \\
&+(\tau^*)^2(t)e^{2\alpha t}\big(|b(\bu_{ht},\bu_h,\bu_{htt})|+ 
|b(\bu_h,\bu_{ht},\bu_{htt})|+\|f_t\|\|\bu_{htt}\|\big).
\end{align}
Here, without loss of generality, we have assumed that $\frac{d}{dt} \tau^*(t) \le 1$. 
{\it To avoid confusion we can always first derive these estimate in the interval $(0,1)$ and then in $(1,t),~t>1$.} For the nonlinear terms, we first note that
\begin{align}\label{dth002}
b(\bu_{ht},\bu_h,\bu_{htt})= \frac{1}{2}(\bu_{ht}\cdot\nabla\bu_h,\bu_{htt})
-\frac{1}{2}(\bu_{ht}\cdot\nabla\bu_{htt},\bu_h).
\end{align}
But we can rewrite the second term as follows: 
(with notations $D_i=\frac{\partial}{\partial x_i}$ and $\bv=(v^1,v^2)$)
\begin{align}\label{dth003}
(\bu_{ht}\cdot\nabla\bu_{htt},\bu_h) &=\sum_{i,j=1}^2\int_{\Omega} u_{ht}^iD_i(u_{htt}^j) u_h^jd\x \\
&= -\sum_{i,j=1}^2\int_{\Omega} D_i(u_{ht}^i)u_{htt}^ju_h^jd\x-\sum_{i,j=1}^2\int_{\Omega} u_{ht}^i u_{htt}^jD_i(u_h^j)d\x \nonumber \\
&= -((\nabla\cdot\bu_{ht})\bu_{htt},\bu_h)-(\bu_{ht}\cdot\nabla\bu_h,\bu_{htt}).  \nonumber
\end{align}
Use (\ref{dth003}) in (\ref{dth002}) and now use (\ref{nonlin1}) to estimate the nonlinear terms.
\begin{align}\label{dth004}
|b(\bu_{ht},\bu_h,\bu_{htt})|+ |b(\bu_h,\bu_{ht},\bu_{htt})| \le \frac{1}{4} 
\|\bu_{htt}\|^2+K\|\bu_{ht}\|_1^2\|\bu_h\|_2^2.
\end{align}
Incorporate (\ref{dth004}) in (\ref{dth001}) and, apply Young's inequality and kickback argument.
Then integrate with respect to time. The resulting double integral is estimated similar to $(3.26)$ of \cite{GP11}. Now using (\ref{dth11})-(\ref{dth12}) we can easily deduce that
\begin{equation*}
 (\tau^*)^2\|\bu_{ht}\|_1^2+\mu e^{-2\alpha t}\int_0^t 
 (\tau^*)^2(s) e^{2\alpha s}\|\bu_{hss}\|^2~ds \le K.
\end{equation*}
We set $\bphi_h=-\tau^*(t) e^{2\alpha t}\td_h^{-1}\bu_{htt}$ in (\ref{dwfjt}). From 
(\ref{nonlin1}) we see that
$$ b(\bu_{ht},\bu_h,\td_h^{-1}\bu_{htt}) \le K\|\bu_{ht}\|^{1/2}\|\bu_{ht}\|_1^{1/2}
   \|\bu_h\|_1\| \bu_{htt}\|_{-1} $$
and therefore
\begin{align*}
 \mu\frac{d}{dt}(\tau^*(t) e^{2\alpha t}\|\bu_{ht}\|^2) +\tau^*(t) e^{2\alpha t} 
\|\bu_{htt}\|_{-1}^2 \le \big(2\alpha \tau^*(t)+1\big) e^{2\alpha t}\|\bu_{ht}\|_1^2 
+C(\mu,\gamma) \tau^*(t) e^{2\alpha t} \|\nabla\bu_h\|^2 +2\|\f_t\|^2 \\
 +C(\mu,\delta)(\int_0^t \beta(t-s) e^{\alpha t}\|\td_h\bu_h(s)\|~ds)^2 
+C(\mu) \tau^*(t) e^{2\alpha t}\Big(\|\nabla\bu_h\|^2\|\bu_{ht}\|^2 
+\|\nabla\bu_{ht}\|^2(1+ \|\bu_h\|\|\nabla\bu_h\|)\Big).
\end{align*}
Integrate with respect to time and multiply by $e^{-2\alpha t}$ to conclude
\begin{equation*}
 \tau^*(t)\|\bu_{ht}\|^2+\mu e^{-2\alpha t}\int_0^t \tau^*(s) e^{2\alpha s} 
\|\bu_{hss}\|_{-1}^2 ds \le K.
\end{equation*}
Finally we set $\bphi_h= -e^{2\alpha t}\td_h^{2}\bu_{htt}$ in (\ref{dwfjt}) and proceed as 
above to arrive at
\begin{equation*}
\|\bu_{ht}\|_{-1}^2+\mu e^{-2\alpha t}\int_0^t e^{2\alpha s} 
\|\bu_{hss}\|_{-2}^2 ds \le K.
\end{equation*}
This completes the rest of the proof.
 \end{proof}
\noindent The following semi-discrete error estimates are proved in \cite{GP11}.
\begin{theorem}\label{errest}
Let $\Omega$ be a convex polygon and let the conditions (${\bf A1}$)-(${\bf A2}$) and (${\bf 
B1}$)-(${\bf B2}$)
be satisfied. Further, let the discrete initial velocity $\bu_{0h}\in \bJ_h$ with
$\bu_{0h}=P_h\bu_0,$ where $\bu_0\in \bJ_1.$ Then,
there exists a positive constant $K$ such that for $0<T<+\infty $ with $t\in (0,T]$
$$ \|(\bu-\bu_h)(t)\|+h\|\nabla(\bu-\bu_h)(t)\|\le Ke^{Kt}h^2t^{-1/2}.$$
Moreover under the assumption of the uniqueness condition, that is,
\begin{equation}\label{uc}
 \frac{N}{\nu^2}\|\f\|_{\infty} < 1~~~\mbox{and}~~N 
=\sup_{\bu, \bv,{\bf w}\in\bH_0^1(\Omega)}\frac{b(\bu,\bv,{\bf 
w})}{\|\nabla\bu\|\|\nabla\bv\| \|\nabla{\bf w}\|},
\end{equation}
where $\nu=\mu+\frac{\gamma}{\delta}$ and $\|\f\|_{\infty} :=\|\f\|_{L^\infty(\R_{+};
 \bL^2(\Omega))}$, we have the following uniform estimate:
$$ \|(\bu-\bu_h)(t)\| \le Kh^2t^{-1/2}. $$
\end{theorem}

\section{Backward Euler Method}
\se
This section deals with fully discrete scheme based on backward Euler method, existence of unique discrete solutions and some uniform in time bounds.
For time discretization, let  $k=\frac{T}{N}$ be the time step with $t_n=nk,~n\ge 0$, $N$ is a positive integer. We define for a sequence $\{\bphi^n\}_{n
\ge 0}\subset\bJ_h,$ the backward difference quotient
$$ \pt\bphi^n=\frac{1}{k}(\bphi^n-\bphi^{n-1}). $$
For any continuous function $\bv(t)$ we set $\bv_n=\bv(t_n)$. 
Since backward Euler method is of first order in time, we choose the right
rectangle rule to approximate the integral term in (\ref{dwfj}).
\begin{equation*}
 q_r^n(\bphi)=k\sum_{j=1}^{n}\beta_{nj}\phi^j\approx \int_0^{t_n}
 \beta(t_n-s)\bphi(s)~ds,
\end{equation*}
where $\beta_{nj}=\beta(t_n-t_j).$ With $w_{nj}= k\beta_{nj}$, it is observed that the 
the right rectangle rule is positive in the sense that
\begin{equation}\label{rrp}
k\sum_{i=1}^n q_r^i(\phi)\phi^i= k\sum_{i=1}^n \sum_{j=0}^i 
\omega_{ij}\phi^j \phi^i \ge 0,~~~~\phi=(\phi^0,\cdots,\phi^N)^T.
\end{equation}
For details, see, McLean and Thom{\'e}e \cite{MT}.
The error incurred due to right rectangle rule in approximating the integral term is
\begin{align}\label{errrr}
\ve_r^n (\phi) & = \int_0^{t_n} \beta(t_n-s)\bphi(s)~ds-k\sum_{j=1}^{n}\beta_{nj}\phi^j \\
&= -\sum_{j=1}^{n}\int_{t_{j-1}}^{t_j} (s-t_{j-1})\frac{\partial}{\partial s}(\beta(t_n-s)
\bphi(s))~ds
\le  k\sum_{j=1}^{n}\int_{t_{j-1}}^{t_j}\Big|\frac{\partial}{\partial s}(\beta(t_n-s) 
\bphi(s))\Big|~ds. \nonumber
\end{align}
We present here a discrete version of \textit{integration by parts} for our subsequent use. For
sequences $\{a_i\}$ and $\{b_i\}$ of real numbers, the following summation by parts holds
\begin{equation}\label{sumbp}
k\sum_{j=1}^i a_jb_j= a_i\widehat{b}_i-k\sum_{j=1}^{i-1} (\pt a_{j+1})\widehat{b}_j, ~~\mbox{ where }
\widehat{b}_i=k\sum_{j=1}^i b_j.
\end{equation}
Here, we observe that for sequences $\{a_i\}$ and $\{b_i\}$ of real numbers, 
\begin{align}\label{chain}
\sigma_i(a_i,\pt b_i) = \pt(\sigma_i(a_i, b_i)) - \sigma_i(\pt a_i, b_i) - (\pt\sigma_i)(a_i, b_i),
\end{align}
where, $\sigma_n=e^{2\alpha t_n}\tau^*(t_n)$ and $\tau^*(t_n)=\min\{1,t_n\}$.
Now, the backward Euler scheme for the semidiscrete Oldroyd problem 
(\ref{dwfh}) is to find $\{\bU^n\}_{1\le n\le N}\in\bH_h$ and $\{P^n\}_{1\le n\le N}\in L_h$ as 
solutions of the recursive nonlinear algebraic equations ($1\le n\le N$):
\begin{equation}\left\{\begin{array}{l}\label{fdbeh}
 (\pt\bU^n,\bphi_h)+\mu a(\bU^n,\bphi_h)+a(q_r^n(\bU),\bphi_h) = (P^n,\nabla\cdot \bphi_h)
 +(\f^n,\bphi_h)-b(\bU^n,\bU^n,\bphi_h),~~~\forall \ \bphi_h\in\bH_h, \\
 (\nabla\cdot\bU^n,\chi_h) =0,~~~\forall \ \chi_h \in L_h,~~~n\ge 0.
\end{array}\right.
\end{equation}
We choose $\bU^0=\bu_{0h}=P_h\bu_0.$ On the other hand, for $\bphi_h\in\bJ_h$ we seek 
$\{\bU^n\}_{1\le n\le N}\in\bJ_h$ such that
\begin{equation}\label{fdbej}
 (\pt\bU^n,\bphi_h)+\mu a(\bU^n,\bphi_h)+a(q_r^n(\bU),\bphi_h)= (\f^n,\bphi_h)
 -b(\bU^n,\bU^n,\bphi_h),~~~\forall\bphi_h\in\bJ_h.
\end{equation}
Using variant of the Brouwer fixed point theorem and standard uniqueness argument, one can show that the discrete problem (\ref{fdbej}) is well-posed. For a proof, we refer to \cite{G11}. Below, we prove
{\it a priori} bounds for the discrete solutions $\{\bU^n\}_{1\le n\le N}.$
We note here that since the bounds proved below are independent of $n,~1\le n\le N$, these bounds are uniform in time, that is, they are still valid as the final time $t_N\to +\infty$.
\begin{lemma}\label{stb}
Let $\alpha_0>0$ be such that for $0<\alpha<\alpha_0$
\begin{equation}\label{alpha}
1+\big(\frac{\mu\lambda_1}{2}\big)k\ge e^{\alpha k}.
\end{equation}
Then, with $0<\alpha<\min\{\alpha_0, \delta,\frac{\mu\lambda_1}{2}\}$ the discrete solution $\bU^n,~1\le n\le N$ of (\ref{fdbej}) satisfies the following estimate:
\begin{align*}
\|\bU^n\|^2+\frac{\mu}{2} e^{-2\alpha t_n}k e^{-\alpha k}\sum_{i=1}^n e^{2\alpha t_i}\|\nabla\bU^i\|^2 \le C
\Big(e^{-2\alpha t_n}\|\bU^0\|^2+\frac{1}{\mu \lambda_1}\|\f\|_{\infty}^2\Big),
\end{align*}
where $\|\f\|_{\infty}=\|\f\|_{L^{\infty}(\R_{+};\bL^2(\Omega))}$.
\end{lemma}
\begin{proof}
Although the proof is similar to \cite[Lemma 9]{PYD06}, we have provided a sketch below for the sake of completeness.
Setting $\tbU^n=e^{\alpha t_n}\bU^n,$ we rewrite (\ref{fdbej}), for $\bphi_h\in\bJ_h,$ as
\begin{equation}\label{tfdbej}
 e^{\alpha t_n}(\pt\bU^n,\bphi_h)+\mu a(\tbU^n,\bphi_h)+e^{-\alpha 
t_n}b(\tbU^n,\tbU^n,
 \bphi_h)+e^{\alpha t_n}a(q_r^n(\bU),\bphi_h)= (\tf^n,\bphi_h).
\end{equation}
Note that
$$ e^{\alpha t_n}\pt\bU^n= e^{\alpha k}\pt\tbU^n-\Big(\frac{e^{\alpha
   k}-1}{k}\Big)\tbU^n. $$
On substituting this in (\ref{tfdbej}) and then multiplying the resulting equation by
$e^{-\alpha k},$ we obtain
\begin{align}\label{stb001}
(\pt\tbU^n,\bphi_h)-\Big(\frac{1-e^{-\alpha k}}{k}\Big)(\tbU^n,\bphi_h)&+e^{-\alpha k}
 \mu a(\tbU^n,\bphi_h)+e^{-\alpha t_{n+1}}b(\tbU^n,\tbU^n,\bphi_h) \nonumber \\
&+\gamma e^{-\alpha k}\sum_{i=1}^n e^{-(\delta-\alpha)(t_n-t_i)} a(\tbU^i,\bphi_h)
 =e^{-\alpha k}(\tf^n,\bphi_h).
\end{align}
Observe that
$$ (\pt\tbU^n,\tbU^n)=\frac{1}{k}(\tbU^n-\tbU^{n-1},\tbU^n) \ge \frac{1}{2k}
   (\|\tbU^n\|^2-\|\tbU^{n-1}\|^2)=\frac{1}{2}\pt\|\tbU^n\|^2 $$
and that the nonlinear term vanishes, that is, $b(\tbU^n,\tbU^n,\tbU^n)=0$. \\
Put $\bphi_h=\tbU^i$ in (\ref{stb001}), for $n=i$ and use Poincar\'e inequality $\|\tbU^i\|^2 \le \frac{1}{\lambda_1}\|\nabla\tbU^i\|^2$ to obtain
\begin{align}\label{stb002}
\frac{1}{2}\pt\|\tbU^i\|^2 -\Big(\frac{1-e^{-\alpha k}}{k\lambda_1}\Big)  \|\nabla\tbU^i\|^2
+e^{-\alpha k}\mu\|\nabla\tbU^i\|^2
+\gamma e^{-\alpha k}k\sum_{j=1}^i e^{-(\delta-\alpha)(t_i-t_j)} a(\tbU^j,\tbU^i)=
e^{-\alpha k}(\tf^i,\tbU^i).
\end{align}
The right-hand side of (\ref{stb002}) can be estimated as
$$ e^{-\alpha k}\|\tf^i\|\|\tbU^i\| \le \frac{1}{4}e^{-\alpha k}\mu \|\nabla\tbU^i\|^2 +\frac{1}{\mu\lambda_1}e^{-\alpha k}\|\tf^i\|^2, $$
so as to obtain from (\ref{stb002})
\begin{align}\label{stb003}
 \pt\|\tbU^i\|^2 + 2e^{-\alpha k}\Big(\frac{3}{4}\mu-\big(\frac{e^{\alpha k}-1}{k\lambda_1}\big)
 \Big)\|\nabla\tbU^i\|^2 
 + 2\gamma e^{-\alpha k}k\sum_{j=1}^i e^{-(\delta-\alpha)(t_i-t_j)} a(\tbU^j,\tbU^i) 
\le \frac{2}{\mu\lambda_1}e^{-\alpha k}\|\tf^i\|^2.
\end{align}
With $0<\alpha < \min\{\alpha_0,\delta,\frac{\mu\lambda_1}{2}\},$ (\ref{alpha}) is satisfied;
which guarantees that $\frac{\mu}{2}\ge \frac{e^{\alpha k}-1}{k\lambda_1}.$ Now multiply (\ref{stb003}) by $k$ and then sum over $i=1$ to $n.$ The resulting double sum is non-negative and hence we obtain
\begin{align} \label{UN-estimate-1}
\|\tbU^n\|^2+ \frac{\mu}{2} k e^{-\alpha k}\sum_{i=1}^n \|\nabla\tbU^i\|^2 \le \|\bU^0\|^2 +\frac{2\|\f\|_{\infty}^2}{\mu\lambda_1}e^{-\alpha k} k\sum_{i=1}^n e^{2\alpha t_i}.
\end{align}
Note that using geometric series, we find, for some $k^*$ in $(0,k)$ that
\begin{equation}\label{e-2alpha}
k\sum_{i=1}^n e^{2\alpha t_i}= e^{2\alpha k}\frac{k}{e^{2\alpha k}-1} e^{2\alpha t_n}= e^{2\alpha(k-k^*)} e^{2\alpha t_n},
\end{equation}
On substituting (\ref{e-2alpha}) in (\ref{UN-estimate-1}),
multiply through out by $e^{-2\alpha t_n}$ to complete the rest of the proof.\\
\end{proof}
\begin{remark}\label{alpha_0}
In Lemma \ref{stb}, such a choice of $\alpha_0>0$ is possible by choosing $\alpha_0 < \frac{\log(1+\frac{\mu\lambda_1}{2}k)}{k}$. Note that for large $k>0$, $\alpha_0$ is small but as $k\to 0$, $\frac{\log(1+\frac{\mu\lambda_1}{2}k)}{k} \to \frac{\mu\lambda_1}{2}$. Thus, with $0<\alpha < \min\{\alpha_0,\delta,\frac{\mu\lambda_1}{2}\}$, the result in Lemma \ref{stb} is valid.
\end{remark}
\noindent
In order to obtain uniform (in time) estimate for the discrete solution $\bU^n$
in Dirichlet norm, we introduce the following notation:
\begin{equation}\label{ubeta}
\bU^n_{\beta}=k\sum_{j=1}^n \beta_{nj}\bU^j, ~n >0;~~\bU^0_{\beta}=0,
\end{equation}
and rewrite (\ref{fdbej}), for $\bphi_h\in\bJ_h,$ as
\begin{equation}\label{fdbej00}
 (\pt\bU^n,\bphi_h)+\mu a(\bU^n,\bphi_h)+b(\bU^n,\bU^n,\bphi_h)
 +a(\bU^n_{\beta},\bphi_h)= (\f^n,\bphi_h).
\end{equation}
Note that
\begin{equation*}
\bU^n_{\beta}= k\gamma\bU^n+e^{-\delta k}\bU^{n-1}_{\beta},
\end{equation*}
and therefore,
\begin{align}\label{une02}
\pt\bU^n_{\beta} =\frac{1}{k}(\bU^n_{\beta}-\bU^{n-1}_{\beta})
=\frac{1}{k}\bU^n_{\beta}-\frac{1}{k}e^{\delta k}(\bU^n_{\beta}-k\gamma\bU^n)
=\gamma e^{\delta k}\bU^n-\frac{(e^{\delta k}-1)}{k}\bU^n_{\beta}.
\end{align}

\begin{lemma}\label{unif.est0}
Under the assumptions of Lemma \ref{stb}, the discrete solution $\bU^n,~1\le n\le N$, of (\ref{fdbej}), satisfies the following uniform estimates:
\begin{equation}\label{unif.est01}
\|\bU^n\|^2+\frac{e^{-\delta k}}{\gamma} \|\nabla\bU^n_{\beta}\|^2 \le e^{-\alpha t_n}\|\bU^0\|^2+ \left(\frac{1-e^{-\alpha t_n}}{\alpha\mu\lambda_1}\right)\|\f\|^2_{\infty} 
=M_{11}^2,
\end{equation}
and
\begin{equation}\label{unif.est02}
k\sum_{n=m}^{m+l} \big(\mu\|\nabla\bU^n\|^2+\frac{\delta}{\gamma} 
\|\nabla\bU^n_{\beta}\|^2  \big) \le M_{11}^2+\frac{l}{\mu\lambda_1} 
\|\f\|^2_{\infty}=M_{12}^2(l),
\end{equation}
where $\bU^n_{\beta}$ is given by (\ref{ubeta}) and $m,l\in\mathbb{N}$.
\end{lemma}

\begin{proof}
In view of (\ref{une02}), we find that
$$ a(\bU^n_{\beta},\bU^n)= \frac{e^{-\delta k}}{\gamma} a(\bU^n_{\beta},\pt\bU^n_{\beta})
+\frac{(1-e^{-\delta k})}{k\gamma}\|\nabla\bU^n_{\beta}\|^2. $$
Now take $\bphi_h=\bU^n$ in (\ref{fdbej00}), for $n=i$, to find
\begin{equation}\label{une03}
\pt\big(\|\bU^i\|^2+\frac{e^{-\delta k}}{\gamma}\|\nabla\bU^i_{\beta}\|^2 \big) 
+2\mu\|\nabla\bU^i\|^2+2\Big(\frac{e^{\delta k}-1}{k}\Big)\frac{e^{-\delta k}}{\gamma}\|\nabla\bU^i_{\beta}\|^2 \le 
 2\|\f^i\|\|\bU^i\| \le    \frac{1}{\mu\lambda_1}\|\f^i\|^2 + \mu\|\nabla\bU^i\|^2.
\end{equation}
A use of the kick back Poincar\'e inequality $ \mu \|\nabla\bU^i\|^2 \ge \mu\lambda_1 \|\bU^i\|^2$ yields
\begin{equation}\label{une04}
\pt\big(\|\bU^i\|^2+\frac{e^{-\delta k}}{\gamma}\|\nabla\bU^i_{\beta}\|^2 \big) 
+ \mu\lambda_1 \|\bU^i\|^2+ 2\Big(\frac{e^{\delta k}-1}{k}\Big) \frac{e^{-\delta k}}{\gamma} 
\|\nabla\bU^i_{\beta}\|^2  \le \frac{1}{\mu\lambda_1}\|\f^i\|^2.
\end{equation}
Multiply the inequality (\ref{une04}) by $e^{\alpha t_{i-1}}$  and 
note that
\begin{align}\label{discrete01}
\pt(e^{\alpha t_i}\bphi^i)= e^{\alpha t_{i-1}}\Big\{\pt\bphi^i +\frac{e^{\alpha k}-1}{k} \bphi^i\Big\}.
\end{align}
Therefore, we obtain from (\ref{une04})
\begin{align*}
\pt \Big(e^{\alpha t_i}\big(\|\bU^i\|^2+\frac{e^{-\delta k}}{\gamma} \|\nabla \bU^i_{\beta}\|^2\big) \Big) 
& + \Big(\mu\lambda_1-\frac{e^{\alpha k}-1}{k} \Big)e^{\alpha t_{i-1}}\|\bU^i\|^2 \\ 
&+ \Big(2\big(\frac{e^{\delta k}-1}{k}\big) -\big(\frac{e^{\alpha k}-1}{k}\big)\Big) \frac{e^{-\delta k}}{\gamma} e^{\alpha t_{i-1}}\|\nabla \bU^i_{\beta}\|^2 
 \le\frac{e^{\alpha t_{i-1}}}{\mu\lambda_1}\|\f\|_{\infty}^2. 
\end{align*}
With $0<\alpha<\min\{\alpha_0, \delta,\mu\lambda_1/2\}$,
the last two terms on the left hand side become non-negative and hence, we drop it. 
Multiply by $k$ and sum over $1$ to $n$ and then multiply the resulting inequality by $e^{-\alpha t_n}.$ Observe that $\bU^0_{\beta}=0$ by definition. This results in the first estimate (\ref{unif.est01}). For the second estimate (\ref{unif.est02}), we multiply (\ref{une03}) by $k,$ sum over $m$ to $m+l$ with $m,l\in\N$ and use (\ref{unif.est01}) to complete the rest of the proof.
\end{proof}

\begin{lemma}\label{unif.est1}
Under the assumptions of Lemma \ref{stb}, the discrete solution $\bU^n,~1\le n\le N$ of (\ref{fdbej}) satisfies the following uniform estimates:
\begin{equation*}
\|\nabla\bU^n\|^2+\frac{e^{-\delta k}}{\gamma}\|\td_h\bU^n_{\beta}\|^2 \le K.
\end{equation*}
\end{lemma}

\begin{proof}
Set $\bphi_h=-\td_h\bU^n$ in (\ref{fdbej00}) and as in the Lemma \ref{unif.est0}, we now 
obtain
\begin{align}\label{une11}
\pt\big(\|\nabla\bU^n\|^2+\frac{e^{-\delta k}}{\gamma}\|\td_h\bU^n_{\beta}\|^2 
\big)+2\mu\|\td_h\bU^n\|^2+2\big(\frac{e^{\delta k}-1}{k}\big)\frac{e^{-\delta k}}{\gamma}\|\td_h\bU^n_{\beta}\|^2 \nonumber\\
 \le 2\|\f^{n}\|\|\td_h\bU^n\|+2|b(\bU^n,\bU^n,-\td_h\bU^n)|.
\end{align}
Use Lemma \ref{nonlin}, \ref{unif.est0} and the Young inequality to estimate the nonlinear term as
\begin{align*}
2|b(\bU^n,\bU^n,-\td_h\bU^n)| 
&\le 2C \|\bU^n\|^{1/2}\|\nabla\bU^n\|^{1/2}\|\nabla\bU^n\|^{1/2}\|\td_h\bU^n\|^{1/2}\|\td_h\bU^n\|\nonumber\\
&\le (\frac{9/2}{\mu})^{3} M_{11}^2\|\nabla\bU^n\|^4 + \frac{\mu}{3} \|\td_h\bU^n\|^{2}
\end{align*} 
Thus we obtain from (\ref{une11})
\begin{align}\label{une12}
 \pt\big(\|\nabla\bU^n\|^2+\frac{e^{-\delta k}}{\gamma}\|\td_h\bU^n_{\beta}\|^2 
\big)+\frac{4\mu} {3}\|\td_h\bU^n\|^2+2\big(\frac{e^{\delta k}-1}{k}\big)\frac{e^{-\delta k}}{\gamma} 
\|\td_h\bU^n_{\beta}\|^2 \nonumber\\
\le \frac{3}{\mu} \|\f\|_{\infty}^2 +(\frac{9/2}{\mu})^{3} M_{11}^2\|\nabla\bU^n\|^4. 
\end{align}
Choose $\gamma_0>0$, 
\begin{equation}\label{une13} 
 \gamma_0\|\nabla\bU^n\|^2 = \gamma_0(\bU^n,-\td_h\bU^n) \le \frac{\mu}{3}\|\td_h\bU^n\|^2 
+\frac{3}{4\mu}\gamma_0^2\|\bU^n\|^2,
\end{equation}
and define
\begin{equation*}
 g^n= \min\Big\{\gamma_0+\mu\lambda_1-(\frac{9}{2\mu})^3 M_{11}^2
 \|\nabla\bU^n\|^2, ~2\big(\frac{e^{\delta k}-1}{k}\big)\Big\},
\end{equation*}
and $E^n :=\|\nabla\bU^n\|^2+\frac{e^{-\delta k}}{\gamma} \|\td_h\bU^n_{\beta}\|^2$ for large enough $\gamma_0$ so that $g^n>0$. Now add the two inequalities (\ref{une12}) and (\ref{une13}) and rewrite the resulting equation as
\begin{equation}\label{une15}
 \pt E^n+g^nE^n \le \frac{3}{\mu}\|\f\|_{\infty}^2+\frac{3}{4\mu}\gamma_0^2\|\bU^n\|^2 
=K_{11}. 
\end{equation}
Let $\{n_{i}\}_{i\in\mathbb{N}}$ and $\{\bar{n}_{i}\}_{i\in\mathbb{N}}$ be two finite subsequences of natural numbers such that
$$ g^{n_{i}}=\gamma_0+\mu\lambda_1-(\frac{9}{2\mu})^3 M_{11}^2\|\nabla\bU^{n_i}\|^2,
~~g^{\bar{n}_{i}}=2\big(\frac{e^{\delta k}-1}{k}\big),~\forall i. $$
If for some $n$,
$$ g^n=\gamma_0+\mu\lambda_1-(\frac{9}{2\mu})^3 M_{11}^2\|\nabla\bU^n\|^2=2\big(\frac{e^{\delta k}-1}{k}\big) $$
then without loss of generality, we assume that $n\in \{\bar{n}_{i}\} $ so as to make
the two subsequence $\{n_{i}\}$ and $\{\bar{n}_{i}\}$ disjoint. Now for $m,l\in\mathbb{N}$, we write
\begin{align}\label{une16}
k\sum_{n=m}^{m+l}g^n = k\sum_{n=m_1}^{m_{l_1}} g^n+k\sum_{n=\bar{m}_1}^{\bar{m}_{l_2}}
g^n =k\sum_{n=m_1}^{m_{l_1}} \left(\gamma_0+\mu \lambda_1-(\frac{9}{2\mu})^3 M_{11}^2  \|\nabla\bU^n\|^2\right)+k\sum_{n=\bar{m}_1}^{\bar{m}_{l_2}} 2\big(\frac{e^{\delta k}-1}{k}\big).
\end{align}
Here, $\{m_1,m_2,\cdots,m_{l_1}\} \subset\{n_{i}\}$ and $\{\bar{m}_1,\bar{m}_2,\cdots,\bar{m}_{l_2}\}
 \subset\{\bar{n}_{i}\}$ such that $l_1+l_2=l+1.$ Note that $l_1$ or $l_2$ could be $0$.
Using Lemma \ref{unif.est0}, we observe that
\begin{align*}
(\frac{9}{2\mu})^3k\sum_{n=m}^{m+l} M_{11}^2 \|\nabla\bU^n\|^2 \le 
\frac{9^3M_{11}^2}{2^3\mu^3}k\sum_{n=m}^{m+l} \|\nabla\bU^n\|^2 \le 
\frac{9^3M_{11}^2}{2^3\mu^4}M_{12}^2(l)=K_{12}(l).
\end{align*}
Therefore, from (\ref{une16}), we find that
\begin{align*}
k\sum_{n=m}^{m+l}g^n \ge kl_1(\gamma_0+\mu \lambda_1)-K_{12}(l_1)+2\big(\frac{e^{\delta k}-1}{k}\big) kl_2.
\end{align*}
We choose $\gamma_0$ such that $kl_1(\gamma_0+\mu \lambda_1)-K_{12}(l_1)=2\big(\frac{e^{\delta k}-1}{k}\big) kl_1$,
assuming $l_1\neq 0$, to arrive at
 \begin{equation}\label{une17}
k\sum_{n=m}^{m+l}g^n \ge 2\big(\frac{e^{\delta k}-1}{k}\big) t_{l+1}.
\end{equation}
By definition of $g^n$, we have equality in (\ref{une17}) and in fact, $g^n=2\big(\frac{e^{\delta k}-1}{k}\big)$.
Now from (\ref{une15}), we obtain
$$ \pt E^n+2\big(\frac{e^{\delta k}-1}{k}\big) E^n \le K_{11}. $$
Multiply the above inequality by $e^{\delta t_{n-1}}$ and as in (\ref{discrete01}), we obtain
$$ \pt(e^{\delta t_n}E^n) + \big(\frac{e^{\delta k}-1}{k}\big)e^{\delta t_{n-1}}E^n
\le K_{11}e^{\delta t_{n-1}}. $$
Multiply by $k$ and sum over $1$ to $n$. Observe that $E^0=\|\nabla\bU^0\|^2$. Finally,
multiply the resulting inequality by $e^{-\delta t_n}$ to find that
$$ E^n \le e^{-\delta t_n}\|\nabla\bU^0\|^2+K. $$
This completes the rest of the proof.
\end{proof}

\begin{remark}
As a consequence of the Lemma 4.3, the following {\it a priori} bound is valid:
\begin{equation*} 
\tau^*(t_n)  \|\td_h\bU^n\|^2 \leq K,~1\le n\le N. 
\end{equation*}
\end{remark}

\section{A Priori Error Estimate}
\se

In this section, we discuss error estimates of the backward Euler method for the Oldroyd 
model (\ref{om})-(\ref{ibc}). For the error analysis, we set, for fixed $n\in\N, ~1\le
n\le N,~\e_n=\bU^n -\bu_h(t_n)=\bU^n-\bu_h^n.$
We now rewrite (\ref{dwfj}) at $t=t_n$ and subtract the resulting one from (\ref{fdbej}) to 
obtain
\begin{align*}
 (\pt\e_n,\bphi_h)+ \mu a(\e_n,\bphi_h)+a(q^n_r(\e),\bphi_h) = E^n (\bu_h)(\bphi_h) 
+\ve_a^n(\bu_h)(\bphi_h)+\Lambda^n_h(\bphi_h),
\end{align*}
where,
\begin{align}
& E^n(\bu_h)(\bphi_h) =  (\bu_{ht}^n,\bphi_h)-(\pt\bu_h^n,\bphi_h) = (\bu_{ht}^n,\bphi_h) 
-\frac{1}{k}\int_{t_{n-1}}^{t_n}(\bu_{hs}, \bphi_h)~ds  
 = \frac{1}{k}\int_{t_{n-1}}^{t_n} (t-t_{n-1})(\bu_{htt},\bphi_h)dt, \label{R1be} \\
& \ve_a^n(\bu_h)(\bphi_h) =  a(\uhb(t_n), \bphi_h)ds -a(q_r^n(\bu_h), \bphi_h)=a(\ve_r^n(\bu_h),\bphi_h), \label{ver} \\
& \Lambda^n_h(\bphi_h) = b(\bu_h^n,\bu_h^n,\bphi_h)-b(\bU^n,\bU^n,\bphi_h) 
 = -b(\bu_h^n,\e_n,\bphi_h)-b(\e_n,\bu_h^n,\bphi_h)-b(\e_n,\e_n,\bphi_h). \label{dLbe}
\end{align}
In order to dissociate the effect of nonlinearity, we first linearized the discrete problem (\ref{fdbej}) and introduce $\{\bV^n\}_{n\ge 1}\in\bJ_h$ as solutions of the following linearized problem:
\begin{equation}\label{fdbejv} 
 (\pt\bV^n,\bphi_h)+\mu a(\bV^n,\bphi_h)+a(q_r^n(\bV),\bphi_h)= (\f^n,\bphi_h)
 -b(\bu_h^n,\bu_h^n,\bphi_h),~~~\forall\bphi_h\in\bJ_h,
\end{equation}
given $\bV^0$ and $\bu_h\in\bJ_h$ as solution of (\ref{dwfj}). It is easy to check the existence and uniqueness of $\{\bV^n\}_{n\ge 1}\in\bJ_h.$ \\
We now split the error as follows:
\begin{align*}
\e_n:=\bU^n-\bu_h^n = (\bU^n-\bV^n)-(\bu_h^n-\bV^n) =: \bta_n-\bxi_n.
\end{align*}
The following equations are satisfied by $\bxi_n$ and $\bta_n,$ respectively:
\begin{align}\label{eebelin}
 (\pt\bxi_n,\bphi_h)+& \mu a(\bxi_n,\bphi_h)+a(q^n_r(\bxi),\bphi_h)
 = -E^n(\bu_h)(\bphi_h)-\ve_a^n(\bu_h)(\bphi_h)
\end{align}
and
\begin{align}\label{eebenl}
 (\pt\bta_n,\bphi_h)+& \mu a(\bta_n,\bphi_h)+a(q^n_r(\bta),\bphi_h) 
 =\Lambda^n_h(\bphi_h).
\end{align}
We first estimate the two terms on the right-hand side of (\ref{eebelin})
in the Lemma below.
\begin{lemma}\label{e-ve}
Let $ r\in \{0,1\}$ and $\alpha$ as defined in  Lemma 4.1. Then
with $E^n$ and $\ve_a^n$ defined, respectively, as (\ref{R1be}) and (\ref{ver}), 
the following estimate holds 
for $n=1,\cdots,N$ and for $\{\bphi_h^i\}_i$ in $\bJ_h$:
\begin{eqnarray}\label{eve1}
 2k\sum_{i=1}^n e^{2\alpha (t_i-t_n)}\Big(E^i(\bu_h)(\bphi_h^i)+\ve_a^i(\bu_h)
(\bphi_h^i)\Big) 
\le Kk^{(1-r)/2}(1+\log\frac{1}{k})^{(1-r)/2}\left({k\sum_{i=1}^n  e^{2\alpha (t_i-t_n)}\|\bphi_h^i\|_{1-r}^2}\right)^{1/2}. \nonumber
\end{eqnarray}
\end{lemma}

\begin{proof}
From (\ref{R1be}), we observe that
\begin{align*}
2k \sum_{i=1}^n  e^{2\alpha (t_i-t_n)} E^i(\bu_h)(\bphi_h^i) 
\le   \left[k^{-1}\sum_{i=1}^n \Big(\int_{t_{i-1}}^{t_i}  e^{\alpha (t_i-t_n)}(t-t_{i-1})\|\bu_{htt}\|_{r-1}\,dt \Big)^2\right]^{1/2}
\left[k\sum_{i=1}^n  e^{2\alpha (t_i-t_n)}\|\bphi_h^i\|_{1-r}^2
\right]^{1/2}. \nonumber
\end{align*}
Now the first term of right hand side of the above inequality can be written as 
\begin{align}\label{eve01a}
& \left[{k^{-1}\sum_{i=1}^n \Big(\int_{t_{i-1}}^{t_i}  e^{\alpha (t_i-t_n)}(t-t_{i-1})\|\bu_{htt}\|_{r-1}dt \Big)^2}\right]^{1/2} \nonumber \\
\le & \left[{k^{-1}e^{-2\alpha t_n}\sum_{i=1}^n \bigg(\int_{t_{i-1}}^{t_i}  \Big((\tau^*)^{-(r+1)/2}(t)(t-t_{i-1})
e^{\alpha (t_i-t)}\Big)  \Big((\tau^*)^{(r+1)/2}(t) e^{\alpha t}\|\bu_{htt}\|_{r-1}\Big)dt}\bigg)^2\right]^{1/2}, 
\end{align}
where, $\tau^*(t)=\min\{1,t\}$. When $r=0$, then the right hand side can be bounded by
\begin{align*}
&\left[{k^{-1}e^{-2\alpha t_n}\sum_{i=1}^n \bigg(\int_{t_{i-1}}^{t_i}  (\tau^*)^{-1}(t)(t-t_{i-1})^2
e^{2\alpha (t_i-t)}dt\bigg)  \bigg(\int_{t_{i-1}}^{t_i} (\tau^*)(t) e^{2\alpha t}\|\bu_{htt}\|_{-1}^2 dt}\bigg)\right]^{1/2} \\
 \le &k^{-1/2}e^{-\alpha t_{n-1}} \left[  \bigg(\int_{0}^{k}  t~ dt\bigg)  \bigg(\int_{0}^{k} (\tau^*)(t) e^{2\alpha t}\|\bu_{htt}\|_{-1}^2 dt\bigg)+k^2\sum_{i=2}^n \bigg(\int_{t_{i-1}}^{t_i}  t^{-1} dt\bigg)  \bigg(\int_{t_{i-1}}^{t_i} (\tau^*)(t) e^{2\alpha t}\|\bu_{htt}\|_{-1}^2 dt \bigg) \right]^{1/2} \\
  \le &K k^{1/2}(1+\log\frac{1}{k})^{1/2}.
\end{align*}
When $r=1$, then, the right hand side is bounded by
\begin{align*}
&\left[{k^{-1}e^{-2\alpha t_n}\sum_{i=1}^n \bigg(\int_{t_{i-1}}^{t_i}  (\tau^*)^{-2}(t)(t-t_{i-1})^2
e^{2\alpha (t_i-t)} (\tau^*)^2(t) e^{2\alpha t}\|\bu_{htt}\|^2 dt}\bigg) \bigg(\int_{t_{i-1}}^{t_i} dt\bigg)\right]^{1/2} \\
 \le &\left[{e^{-2\alpha t_{n}}e^{2\alpha k}\sum_{i=1}^n \bigg(\int_{t_{i-1}}^{t_i}  (\tau^*)^2(t) e^{2\alpha t}\|\bu_{htt}\|^2 dt}\bigg)\right]^{1/2}   \le K.
\end{align*}
This completes the proof of the first half. For the remaining part, we observe from (\ref{ver}) and
(\ref{errrr}) that
\begin{align}\label{eve02}
& 2k\sum_{i=1}^n  e^{2\alpha (t_i-t_n)} \ve_a^i(\bu_h)(\bphi_h^i) \le \left[{k\sum_{i=1}^n  e^{2\alpha (t_i-t_n)}\|\bphi_h^i\|_{1-r}^2}\right]^{1/2}~\times \\
& \left[{4k\sum_{i=1}^n \Big(\sum_{j=1}^i \int_{t_{j-1}}^{t^j} e^{\alpha (t_i-t_n)}(t-t_{j-1})
\beta(t_i-t)\{\delta\|\bu_h\|_{r+1}+\|\bu_{ht}\|_{r+1}
\}\,dt\Big)^2}\right]^{1/2}. \nonumber 
\end{align}
In Lemma \ref{dth2}, we find that the estimates of $\|\bu_{htt}\|_{r-1}$ and $\|\bu_{ht}\|_{r+1}$ are similar, in fact, the powers of $t_i$ are same. Therefore,the right-hand side of (\ref{eve02}) involving $\|\bu_{ht}\|_{r+1}$ can be estimated
similarly as in (\ref{eve01a}). The terms involving $\|\bu_h\|_{r+1}$ are clearly easy to estimate. But for the sake of completeness, we provide  the case, when $ r=0$ as
 \begin{align*}
& 4\delta^2 k\sum_{i=1}^n \Big(\sum_{j=1}^i \int_{t_{j-1}}^{t^j} e^{\alpha (t_i-t_n)}
(t-t_{j-1})\beta(t_i-t) \|\nabla\bu_h\|~dt\Big)^2 \nonumber \\
\le & ~4\gamma^2\delta^2 e^{-2\alpha t_n} k^3\sum_{i=1}^n e^{-2(\delta-\alpha)t_i} \Big(\sum_{j=1}^i \int_{t_{j-1}}^{t^j}e^{(\delta-\alpha)t}\|\nabla\tilde{\bu}_h\|~dt\Big)^2 \nonumber \\
 \le & ~4\gamma^2\delta^2 e^{-2\alpha t_n} k^3\sum_{i=1}^n e^{-2(\delta-\alpha)t_i} \Big(\int_0^{t_i} e^{2(\delta-\alpha)s}ds\Big) \Big(\int_0^{t_i} \|\nabla\tilde{\bu}_h(s)\|^2 ds\Big) \nonumber \\
 \le & \frac{2\gamma^2\delta^2}{2 (\delta-\alpha)} e^{-2\alpha t_n}k^3 \sum_{i=1}^n e^{2(\delta-\alpha)k}\big(Ke^{2\alpha t_i}\big) \le Kk^3 e^{2\delta k}.
 \end{align*}
This completes the rest of the proof.
\end{proof}
 
\begin{lemma}\label{pree}
 Assume (${\bf A1}$)-(${\bf A2}$) and a space discretization scheme that satisfies conditions (${\bf B1}$)-(${\bf B2}$). 
Let $\alpha_0>0$ be such that $0<\alpha<\min\{\alpha_0, \delta,\frac{\mu\lambda_1}{2}\}$, (\ref{alpha}) be satisfied.
Further, assume that $\bu_h(t)$ and $\bV^n$ satisfy  (\ref{dwfj}) and (\ref{fdbejv}), respectively. 
Then, there is a positive
constant $K$ such that, $\bxi_n=\bu_h^n-\bV^n,~1\le n\le N$, satisfy
\begin{eqnarray} \label{pree1}
 \|\bxi_n\|^2&+& e^{-2\alpha t_n} k\sum_{i=1}^n e^{2\alpha t_i} \|\bxi_i\|_1^2 \le  Kk\big(1+\log\frac{1}{k}\big), \\
 \|\bxi_n\|_1^2&+& e^{-2\alpha t_n} k\sum_{i=1}^n e^{2\alpha t_n} \{\|\bxi_i\|_2^2+ \|\pt\bxi_i\|^2\} \le K.  \label{pree2}
\end{eqnarray}
\end{lemma}

\begin{proof}
For $n=i,$ we put $\bphi_h=\bxi_i$ in (\ref{eebelin}) and with the observation
$$ (\pt\bxi_i,\bxi_i)=\frac{1}{k}(\bxi_i-\bxi_{i-1},\bxi_i) \ge \frac{1}{2k}(\|\bxi_i\|^2
   -\|\bxi_{i-1}\|^2)=\frac{1}{2}\pt\|\bxi_i\|^2, $$
we find that
\begin{equation}\label{pree01}
 \pt\|\bxi_i\|^2+2\mu\|\nabla\bxi_i\|^2+2a(q^i_r(\bxi),\bxi_i) \le -2E^i(\bu_h)(\bxi_i) -2\ve_a^i(\bu_h)(\bxi_i).
\end{equation}
Multiply (\ref{pree01}) by $ke^{2\alpha t_i}$ and sum over $1\le i\le n\le N$ to obtain
\begin{align}\label{pree02}
\|\tbxi_n\|^2 - 2\sum_{i=1}^{n-1} (e^{2\alpha k}-1)\|\tbxi_i\|^2 +2\mu k\sum_{i=1}^n \|\nabla\tbxi_i\|^2 & 
\le  - 2k\sum_{i=1}^n e^{2\alpha t_i}\Big\{E^i(\bu_h)(\bxi_i) +\ve_a^i(\bu_h)(\bxi_i)\Big\} \nonumber \\
& \le  \frac{\mu}{2} k\sum_{i=1}^n \|\nabla\tbxi_i\|^2+Kk \big(1+\log\frac{1}{k}\big)
 e^{2\alpha t_{n+1}}.
\end{align}
Recall that $\tilde{v}(t)=e^{\alpha t}v(t)$.
Note that we have dropped the quadrature term on the left hand-side of (\ref{pree01}) after summation as it is non-negative. Finally, we have used Lemma \ref{e-ve} for $r=0$. Now, we use poincar\'e inequality in the second tern on the left hand side of (\ref{pree02}) to obtain
\begin{align}\label{pree06}
\|\tbxi_n\|^2 &+ 2\Big(\frac{3}{4}\mu - \big(\frac{e^{2\alpha k}-1}{k\lambda_1}\big)\Big)  k\sum_{i=1}^n 
\|\nabla\tbxi_i\|^2\le Kk \big(1+\log\frac{1}{k}\big)e^{2\alpha t_{n+1}}.
\end{align}
Similar to Lemma \ref{stb}, one can show that $2\Big(\frac{3}{4}\mu - \big(\frac{e^{2\alpha k}-1}{k\lambda_1}\big)\Big)\ge\frac{\mu}{2}>0$.
Now multiply (\ref{pree06}) by $e^{-2\alpha t_n}$ to establish (\ref{pree1}).
Next, for $n=i,$ we put $\bphi_h=-\td_h\bxi_i$ in (\ref{eebelin}) and follow
as above to obtain the first part of (\ref{pree2}), that is,
$$ \|\bxi_n\|_1^2+e^{-2\alpha t_n}k\sum_{i=1}^n e^{2\alpha t_i}\|\bxi_i\|_2^2 \le K. $$
Here, we have used (\ref{eve1}) for $r=1$  replacing
$\bphi_h^{i}$ by $\td_h\bxi_i.$
Finally, for $n=i,$ we put $\bphi_h=\pt\bxi_i$ in (\ref{eebelin}) to find that
\begin{align}\label{pree07}
2\|\pt\bxi_i\|^2+\mu\pt\|\bxi_i\|_1^2 \le -2a(q_r^i(\bxi),\pt\bxi_i)
-2E^i(\bu_h)(\pt\bxi_i) -2\ve_a^i(\bu_h)(\pt\bxi_i).
\end{align}
Multiply (\ref{pree07}) by $ke^{2\alpha t_i}$ and sum over $1\le i\le n\le N$.
As has been done earlier, using (\ref{eve1}) for $r=0$, we obtain
\begin{align}\label{pree07a}
\sum_{i=1}^n ke^{2\alpha t_i} \big\{2E^i(\bu_h)(\pt\bxi_i)+2\ve_a^i(\bu_h)(\pt\bxi_i) \big\}
\le \frac{k}{2} \sum_{i=1}^n e^{2\alpha t_i}\|\pt\bxi_i\|^2+K.
\end{align}
The only difference is that the resulting double sum (the term involving $q_r^i$)
is no longer non-negative and hence, we need to estimate it. Note that
\begin{align}\label{pree08}
2k\sum_{i=1}^n e^{2\alpha t_i} a(q_r^i(\bxi),\pt\bxi_i)=2\gamma k^2\sum_{i=1}^n
\sum_{j=1}^i e^{-(\delta-\alpha)(t_i-t_j)} a(\tbxi_j,e^{\alpha t_i}\pt\bxi_i) \\
\le \frac{k}{2} \sum_{i=1}^n e^{2\alpha t_i}\|\pt\bxi_i\|^2
+K(\gamma) k\sum_{i=1}^n \Big(k\sum_{j=1}^i e^{-(\delta-\alpha)(t_i-t_j)} 
\|\td_h\tbxi_j\|\Big)^2. \nonumber
\end{align}
Using change of variable and change of order of double sum, we obtain
\begin{align}\label{pree09}
I & := K(\gamma) k\sum_{i=1}^n \Big(k\sum_{j=1}^i e^{-(\delta-\alpha)(t_i-t_j)} 
\|\td_h\tbxi_j\|\Big)^2 \nonumber\\
& \le K(\alpha,\gamma) e^{(\delta-\alpha)k} k^2 \sum_{i=1}^n k\sum_{j=1}^i e^{-(\delta-\alpha)(t_i-t_j)}\|\td_h\tbxi_j\|^2\nonumber\\
& \le K(\alpha,\gamma) e^{(\delta-\alpha)k}k^2\sum_{i=1}^n k\sum_{l=1}^{i} e^{-(\delta-\alpha)t_{l-1}}\|\td_h\tbxi_{i-l+1}\|^2~~~\mbox{for} ~l=i-j \nonumber\\
& \le K(\alpha,\gamma) e^{(\delta-\alpha)k}k \Big(k\sum_{l=1}^{n-1} e^{-(\delta-\alpha)t_l}\Big)\Big(k\sum_{j=1}^n \|\td_h\tbxi_j\|^2\Big) \le K. 
\end{align}
%
%
Combining (\ref{pree08})-(\ref{pree09}), we find that
\begin{align}\label{pree10}
2k\sum_{i=1}^n e^{2\alpha t_i} a(q_r^i(\bxi),\pt\bxi_i)\le \frac{k}{2}
\sum_{i=1}^n e^{2\alpha t_i}\|\pt\bxi_i\|^2+K.
\end{align}
Incorporating (\ref{pree07a}) and (\ref{pree10}), we obtain
\begin{align*}
k\sum_{i=1}^n e^{2\alpha t_i}\|\pt\bxi_i\|^2+\mu\|\tbxi_n\|_1^2 \le K+\mu
k \sum_{i=1}^{n-1} \frac{ (e^{2\alpha k}-1)}{k}\|\tbxi_i\|_1^2.
\end{align*}
Use (\ref{pree1}) and the fact that $(e^{2\alpha k}-1)/k \le K(\alpha)$ to
complete the rest of the proof.
\end{proof}


\noindent
Lemma \ref{pree} provides us with the estimate $\|\bxi_i\| \le Kk^{1/2}\big(1+\log
\frac{1}{k}\big)^{1/2}$, which is suboptimal. Before attempting to improve it, we need a few intermediate estimates. First, we present optimal $l^2(\bL^2)$-norm estimate of $\bxi_i$.
Analogous to the semi-discrete case, we resort to duality argument to obtain the same. \\
Consider the following backward problem: For a given
$\bW_n$ and $\g_i,$ let $\bW_i,~n\ge i\ge 1 $ satisfy 
\begin{equation}\label{bwprob}
(\bphi_h,\pt\bW_i)-\mu a(\bphi_h, \bW_i)-k\sum_{j=i}^n \beta(t_j-t_i) a(\bphi_h, \bW_j) 
=(\bphi_h, e^{2\alpha t_i} \g_i),~~\forall\bphi_h\in \bJ_h.
\end{equation}
By denoting $\bar{\bW}_i= \bW_{n-i}$, we shall obtain a forward problem in $\{\bar{\bW}_i\}$, similar to (\ref{dwfj}), but linear. Following the line of argument used to prove Lemma \ref{stb}, the following {\it a priori} estimates are easy to derive.

\begin{lemma}\label{bwest}
Let the assumptions (${\bf A1}$), (${\bf A2}$), (${\bf B1}$) and (${\bf B2}$) hold. Then,
the following estimates hold under appropriate assumptions on $\bW_n$ and $g$:
\begin{equation*}
\|\bW_0\|_r^2+ k\sum_{i=1}^n e^{-2\alpha t_i}\{\|\bW_i\|_{r+1}+ \|\pt\bW_i\|_{r-1}\} \le K\big\{\|\bW_n\|_r^2+k\sum_{i=1}^n e^{2\alpha t_i} \|\g_i\|_{r-1}^2\big\},
\end{equation*}
where $r\in \{0,1\}$.
\end{lemma}

\begin{lemma}\label{neg}
Under the assumptions of Lemma \ref{bwest}, the following estimate holds:
\begin{equation*}
e^{-2\alpha t_n} k \sum_{i=1}^n e^{2\alpha t_i} \|\bxi_i\|^2 \le Kk^2. 
\end{equation*}
\end{lemma}

\begin{proof}
With 
$ \bW_n=(-\td_h)^{-1}\bxi_n,~~\g_i=\bxi_i~\forall i,$ 
we choose $\bphi_h=\bxi_i$ in (\ref{bwprob}) and use (\ref{eebelin}) to obtain
\begin{align}\label{neg01}
e^{2\alpha t_i}\|\bxi_i\|^2 
&=(\bxi_i,\pt\bW_i)-\mu a(\bxi_i, \bW_i)-k\sum_{j=i}^n \beta(t_j-t_i) a(\bxi_i, \bW_j) \nonumber \\
&= \pt (\bxi_i,\bW_i)-(\pt\bxi_i,\bW_{i-1})-\mu a(\bxi_i, \bW_i)-k\sum_{j=i}^n \beta(t_j-t_i) a(\bxi_i, \bW_j)\nonumber \\
&= \pt (\bxi_i,\bW_i)+k(\pt\bxi_i,\pt\bW_i)+k\sum_{j=1}^i \beta(t_i-t_j) a(\bxi_j, \bW_i)+E^i(\bu_h)(\bW_i) \nonumber\\
&~~~~+\ve_a^i(\bu_h)(\bW_i)-k\sum_{j=i}^n \beta(t_j-t_i) a(\bxi_i, \bW_j).  
\end{align}
Multiply (\ref{neg01}) by $k$ and sum over $1\le i\le n$. Observe that the resulting two 
double sums cancel out (change of order of double sum). Therefore, we find that
\begin{align}\label{neg02}
k\sum_{i=1}^n e^{2\alpha t_i}\|\bxi_i\|^2+\|\bxi_n\|_{-1}^2 =k\sum_{i=1}^n \big[ 
k(\pt\bxi_i,\pt\bW_i)+E^i(\bu_h)(\bW_i)+\ve_a^i(\bu_h) (\bW_i)\big].
\end{align}
From (\ref{R1be}), we observe that
\begin{align}\label{neg03}
k\sum_{i=1}^n E^i(\bu_h)(\bW_i) \le k\sum_{i=1}^n \frac{1}{k}\int_{t_{i-1}}^{t_i} 
(s-t_{i-1})\|\bu_{hss}\|_{-2}\|\bW_i\|_2 \nonumber \\
\le k e^{\alpha k}\Big(\int_0^{t_n} e^{2\alpha s}\|\bu_{hss}\|_{-2}^2 
ds\Big)^{1/2}\Big(k\sum_{i=1}^n e^{-2\alpha t_i}\|\bW_i\|_2^2\Big)^{1/2}.
\end{align}
Similar to (\ref{pree06}), we obtain
\begin{align}\label{neg04}
k\sum_{i=1}^n \ve_a^i(\bu_h)(\bW_i)\le K\Big(k^3\sum_{i=1}^n \int_0^{t_i} e^{2\alpha s} 
(\|\bu_h\|^2+\|\bu_{hs}\|^2)~ds\Big)^{1/2}\Big(k\sum_{i=1}^n e^{-2\alpha t_i} 
\|\bW_i\|_2^2\Big)^{1/2},
\end{align}
and
\begin{align}\label{neg05}
k\sum_{i=1}^n k(\pt\bxi_i,\pt\bW_i)\le k\Big(k\sum_{i=1}^n e^{2\alpha t_i} \|\pt \bxi_i\|^2\Big)^{1/2}\Big(k\sum_{i=1}^n e^{-2\alpha t_i}\|\pt\bW_i\|^2\Big)^{1/2}.
\end{align}
Incorporating (\ref{neg03})-(\ref{neg05}) in (\ref{neg02}), and using Lemmas \ref{dth2}, \ref{pree} and \ref{bwest}, we find that
\begin{align*}
k\sum_{i=1}^n e^{2\alpha t_i}\|\bxi_i\|^2+\|\bxi_n\|_{-1}^2 \le Kk^2 e^{2\alpha t_n}.
\end{align*}
This concludes the rest of the proof.
\end{proof}

\noindent We recall that the estimate of $\|\bxi_n\|$ given in Lemma \ref{pree} is suboptimal, due to non-smooth initial data. The standard technique to improve it is to use a discrete weight function $t_n$. But this would not go through in our case, due to the quadrature term, without some additional tool. This will be clear in our next result.
\begin{lemma}\label{l2eebxi}
 Under the assumptions of Lemma \ref{pree}, the following holds:
\begin{equation}\label{l2eebxi1}
t_n \|\bxi_n\|^2 +e^{-2\alpha t_n} k\sum_{i=1}^n \sigma_i\|\pt\bxi_i\|_{-1}^2 \le 
Kk^2(1+\log \frac{1}{k}),~~1\le n\le N,
\end{equation}
where $\sigma_i=t_i e^{2\alpha t_i}$.
\end{lemma}

\begin{proof}
 From (\ref{eebelin}) with $n=i$ and $\bphi_h=\sigma_i(-\td_h)^{-1}\pt\bxi_i$, we obtain
\begin{align}\label{teebelin} 
2\sigma_i\|\pt\bxi_i\|_{-1}^2 &+\mu \sigma_i\pt\|\bxi_i\|^2 +2\sigma_i a(q_{r}^i(\bxi), (-\td_h)^{-1}\pt\bxi_i) \nonumber\\
    &\le -2E^i(\bu_h)(\sigma_i(-\td_h)^{-1}\pt\bxi_i)- 2\ve_{a}^i(\bu_h) (\sigma_i (-\td_h)^{-1}\pt\bxi_i).
\end{align}
We multiply (\ref{teebelin}) by $k$ and sum it over $1$ to $n$ and using the fact 
$$k\sum_{i=1}^n\sigma_i\pt\|\bxi_i\|^2 \geq \sigma_n\|\bxi_n\|^2 - k\sum_{i=1}^{n-1} e^{2\alpha t_{i}}\|\bxi_{i}\|^2$$
to find that
\begin{align}\label{l2eebxi01}
\mu\sigma_n\|\bxi_n\|^2 + 2 k\sum_{i=1}^n 
\sigma_i\|\pt\bxi_i\|_{-1}^2+2k\sum_{i=1}^n \sigma_i a(q_{r}^i(\bxi), (-\td_h)^{-1}\pt\bxi_i)
\le  k\sum_{i=1}^{n-1}e^{2\alpha t_i}\|\bxi_i\|^2 \nonumber \\
-2 k\sum_{i=1}^n E^i(\bu_h)(\sigma_i(-\td_h)^{-1}\pt\bxi_i)-2k\sum_{i=1}^n \ve_{a}^i(\bu_h) (\sigma_i (-\td_h)^{-1}\pt\bxi_i).
\end{align}
A use of Cauchy-Schwarz's inequality with (\ref{R1be}) yields
\begin{align}\label{l2eebxi02}
k\sum_{i=1}^{n} \sigma_i E^i(\bu_h)((-\td_h)^{-1}\pt\bxi_i) &\le k\sum_{i=1}^n  \sigma_i \frac{1}{k}\int_{t_{i-1}}^{t_i} (s-t_{i-1})\|\bu_{hss}\|_{-1} ds \|\pt\bxi_i\|_{-1} \nonumber\\
	&\le K  k\sum_{i=1}^n   \sigma_i \bigg(\int_{t_{i-1}}^{t_i} \|\bu_{hss}\|_{-1} ds\bigg)^2 + \epsilon k\sum_{i=1}^n \sigma_i \|\pt\bxi_i\|_{-1}^2 \nonumber\\
	&\le K  k^2\sum_{i=1}^n    \int_{t_{i-1}}^{t_i} \sigma_i\|\bu_{hss}\|_{-1}^2 ds + \epsilon k\sum_{i=1}^n \sigma_i \|\pt\bxi_i\|_{-1}^2 \nonumber\\
	&\le K k^2 e^{2\alpha t_n} + \epsilon k\sum_{i=1}^n \sigma_i \|\pt\bxi_i\|_{-1}^2.
\end{align}
Now from (\ref{ver}) and (\ref{chain}), it follows that
\begin{align}
k\sum_{i=1}^{n} \sigma_i \ve_{a}^i(\bu_h)((-\td_h)^{-1}\pt\bxi_i) &= k\sum_{i=1}^{n} \sigma_i a(\ve_{r}^i(\bu_h),(-\td_h)^{-1}\pt\bxi_i)= k\sum_{i=1}^{n} \sigma_i (\ve_{r}^i(\bu_h),\pt\bxi_i) \nonumber \\
	&=\sigma_n (\ve_{r}^n(\bu_h),\bxi_n) - k\sum_{i=1}^{n} (\pt\sigma_i) (\ve_{r}^i(\bu_h),\bxi_i) - k\sum_{i=1}^{n} \sigma_i (\pt \ve_{r}^i(\bu_h),\bxi_i).
\end{align}
Using (\ref{errrr}) and Cauchy-Schwarz's inequality, we can bound
\begin{align}
\sigma_n (\ve_{r}^n(\bu_h),\bxi_n) &\le K k^2\sigma_n\Big(\sum_{i=1}^n \int_{t_{i-1}}^{t_i} \beta(t_n-s)\{\delta\|\bu_h\|+\|\bu_{hs}\|\} ds \Big)^2 + \epsilon \sigma_n \|\bxi_n\|^2 
   \le K k^2 e^{2\alpha t_n} + \epsilon \sigma_n \|\bxi_n\|^2.
\end{align}
A use of (\ref{errrr}) with Lemma \ref{neg} shows
\begin{align}
k\sum_{i=1}^{n} (\pt\sigma_i) (\ve_{r}^i(\bu_h),\bxi_i) &\le k\sum_{i=1}^{n} e^{2\alpha t_i} (\ve_{r}^i(\bu_h),\bxi_i)) \nonumber\\
	&\le K k^2 k\sum_{i=1}^{n} e^{2\alpha t_i} \int_{0}^{t_i} \{\|\bu_h\|^2+\|\bu_{hs}\|^2\} ds + K k\sum_{i=1}^{n} e^{2\alpha t_i} \|\bxi_i\|^2 
	 \le K k^2 e^{2\alpha t_n},
\end{align}
and
\begin{align}
k\sum_{i=1}^{n} \sigma_i (\pt \ve_{r}^i(\bu_h),\bxi_i) 	&\le K k^2 k\sum_{i=1}^{n} e^{2\alpha t_i} \int_{t_{i-1}}^{t_i} \{\|\bu_h\|^2+\|\bu_{hs}\|^2\} ds + K k\sum_{i=1}^{n} e^{2\alpha t_i} \|\bxi_i\|^2 
	 \le K k^2 e^{2\alpha t_n}.\label{l2eebxi02a}
\end{align}
Incorporate (\ref{l2eebxi02})-(\ref{l2eebxi02a}) in (\ref{l2eebxi01}) and use Lemma \ref{neg} to arrive at
\begin{align}\label{l2eebxi03}
\sigma_n\|\bxi_n\|^2 + k\sum_{i=1}^n 
\sigma_i\|\pt\bxi_i\|_{-1}^2+ 2k\sum_{i=1}^n \sigma_i a(q_{r}^i(\bxi), (\td_h)^{-1}\pt\bxi_i)
\le Kk^2  e^{2\alpha t_n}.
\end{align}
Unfortunately, the last term on the left-hand side of (\ref{l2eebxi03}) is no longer non-negative due to the presence of $t_i$. We need an estimate for it and standard tools deployed till now will only result in suboptimal estimate. \\
With the assumption of the estimate
\begin{align}\label{quad.hat}
2k\sum_{i=1}^n \sigma_i a(q_{r}^i(\bxi), (\td_h)^{-1}\pt\bxi_i) \le \ve k\sum_{i=1}^n \sigma_i 
\|\pt\bxi_i\|_{-1}^2+K k\sum_{i=1}^n e^{2\alpha t_i}\|\nabla\hbxi_i\|^2
\end{align}
and with appropriate choice of $\ve >0$, we complete the rest of the proof.
\end{proof}
\noindent
We now have a task at hand of proving (\ref{quad.hat}) and for that purpose, we introduce the
`hat operator', that is,
\begin{equation}\label{sum0}
{\widehat{\bphi}}_h^n := k\sum_{i=1}^n \bphi_h^i.
\end{equation}
This can be considered as discrete integral operator. We first observe, using (\ref{sumbp}), 
that
\begin{align*}
 k\sum_{j=1}^i \beta(t_i-t_j) \bphi_j = \gamma e^{-\delta t_i}k\sum_{j=1}^i e^{\delta t_j} 
\bphi_j 
=  \gamma e^{-\delta t_i} \Big\{e^{\delta t_i}\widehat{\bphi}_i- k\sum_{j=1}^{i-1} 
(\frac{e^{\delta t_{j+1}}-e^{\delta t_j}}{k}) \widehat{\bphi}_j \Big\}= \pt^i \Big\{k 
\sum_{j=1}^i \beta(t_i-t_j) \widehat{\bphi}_j \Big\}.
\end{align*}
Here, $\pt^i$ means the backward difference formula with respect to $i$. Now rewrite the equations 
(\ref{eebelin}) (for $n=i$) as follows: 
\begin{align}\label{eebelin1}
 (\pt\bxi_i,\bphi_h)+& \mu a(\bxi_i,\bphi_h)+\pt^i \Big\{k\sum_{j=1}^i \beta(t_i-t_j)  a 
(\hbxi_j,\bphi_h) \Big\} = -E^i(\bu_h)(\bphi_h)-\ve_a^i(\bu_h)(\bphi_h).
\end{align}
We multiply (\ref{eebelin1}) by $k$ and sum over $1$ to $n$.  Using the fact that 
$\pt\hbxi_n=\bxi_n,$ we observe that
\begin{align}\label{eebelin1i}
 (\pt\hbxi_n,\bphi_h)+ \mu a(\hbxi_n,\bphi_h)+ a(q^n_r(\hbxi),\bphi_h) = -k\sum_{i=1}^n 
\big(E^i(\bu_h)(\bphi_h)+\ve_a^i(\bu_h)(\bphi_h)\big).
\end{align}

\begin{lemma}\label{ieens}
 Under the assumptions of Lemma \ref{pree}, the following estimate holds:
\begin{align}\label{ieens1}
\|\hbxi_n\|^2+ e^{-2\alpha t_n} k\sum_{i=1}^n e^{2\alpha t_i}\|\nabla\hbxi_i\|^2 \le  Kk^2(1+\log\frac{1}{k}),~1\le n\le N.
\end{align}
\end{lemma}

\begin{proof}
Choose $\bphi_h=\hbxi_i$ in (\ref{eebelin1i}) for $n=i$, multiply by $k e^{2\alpha t_i}$ and 
then sum over $1\le i\le n$. We drop the third term on the left hand-side of the resulting 
inequality due to non-negativity.
\begin{align}\label{ieens01}
e^{2\alpha t_n}\|\hbxi_n\|^2+\mu k\sum_{i=1}^n e^{2\alpha t_i}\|\nabla\hbxi_i\|^2 \le 
k\sum_{i=1}^n e^{2\alpha t_i} k\sum_{j=1}^i \big(|E^j(\bu_h)(\hbxi_i)| 
+|\ve_a^j(\bu_h)(\hbxi_i)|\big).
\end{align}
From (\ref{R1be}) and use the similar technique of Lemma \ref{e-ve} to find that
$$ k\sum_{j=1}^i |E^j(\bu_h)(\hbxi_i)| \le \big(\sum_{j=1}^i \int_{t_{j-1}}^{t_j} (s-t_{j-1}) 
\|\bu_{hss}\|_{-1}ds \big)\|\nabla\hbxi_i\| \le Kk(1+\log\frac{1}{k})^{1/2} e^{-\alpha k} \|\nabla\hbxi_i\|. $$
Therefore
\begin{align}\label{ieens03}
k\sum_{i=1}^n e^{2\alpha t_i} k\sum_{j=1}^i |E^j(\bu_h)(\hbxi_i)| \le \frac{\mu}{4} k
 \sum_{i=1}^n e^{2\alpha t_i}\|\nabla\hbxi_i\|^2+Kk^2(1+\log\frac{1}{k}) e^{2\alpha t_n}.
\end{align}
Similarly
\begin{align}\label{ieens04b}
k\sum_{i=1}^n e^{2\alpha t_i} k\sum_{j=1}^i |\ve^j_a(\bu_h)(\hbxi_i)| \le
 \frac{\mu}{4} k \sum_{i =1}^n e^{2\alpha t_i}\|\nabla\hbxi_i\|^2+Kk^2
 (1+\log\frac{1}{k}) e^{2\alpha t_n}.
\end{align}
Incorporate (\ref{ieens03})-(\ref{ieens04b}) in (\ref{ieens01}) to complete 
the rest of the proof.
\end{proof}
\noindent
We are now in a position to prove (\ref{quad.hat}). We rewrite the quadrature term in terms of `hat' and take advantage of the optimal result (\ref{ieens1}).
\begin{lemma}\label{quad.hat.est}
Under the assumption of Lemma \ref{pree}, the following estimate holds:
$$ 2k\sum_{i=1}^n \sigma_i a(q_{r}^i(\bxi), (\td_h)^{-1}\pt\bxi_i) \le \ve k\sum_{i=1}^n \sigma_i \|\pt\bxi_i\|_{-1}^2+Kk^2 (1+\log\frac{1}{k}) e^{2\alpha t_n}. $$
\end{lemma}
\begin{proof}
Using (\ref{sumbp}), we first rewrite the quadrature term as
\begin{align}\label{qhe01}
2k\sum_{i=1}^n \sigma_i a(q_{r}^i(\bxi), (\td_h)^{-1}\pt\bxi_i)= 2k\sum_{i=1}^n \gamma a(\hbxi_i, 
\sigma_i(\td_h)^{-1}\pt\bxi_i)-2k\sum_{i=2}^n k\sum_{j=1}^{i-1} \pt \beta(t_i-t_j) a(\hbxi_j, 
\sigma_i(\td_h)^{-1}\pt\bxi_i).
\end{align}
The first term can be handled as follows (for some $\ve >0$):
\begin{align}\label{qhe02}
2k\sum_{i=1}^n \gamma a(\hbxi_i, \sigma_i(\td_h)^{-1}\pt\bxi_i) \le \ve k\sum_{i=1}^n \sigma_i 
\|\pt\bxi_i\|_{-1}^2+c(\ve, \mu,\gamma) k\sum_{i=1}^n e^{2\alpha t_i} 
\|\nabla\hbxi_i\|^2.
\end{align}
For the second term, using similar technique as in (\ref{pree09}), we observe that
\begin{align}\label{qhe03}
 2k\sum_{i=2}^n k\sum_{j=1}^{i-1} \pt \beta(t_i-t_j) a(\hbxi_j, \sigma_i(\td_h)^{-1}\pt\bxi_i)
&\le \ve k\sum_{i=1}^n \sigma_i \|\pt\bxi_i\|_{-1}^2 + Kk\sum_{i=2}^n \Big(k\sum_{j=1}^{i-1} e^{-\delta(t_i-t_j)} \big(\frac{e^{\delta k-1}}{k}\big) e^{\alpha t_i}\|\nabla\hbxi_j\|\Big)^2 \nonumber\\
&\le \ve k\sum_{i=1}^n \sigma_i \|\pt\bxi_i\|_{-1}^2+Kk\sum_{i=1}^n e^{2\alpha t_i}\|\nabla\hbxi_i\|^2. 
\end{align}
With (\ref{qhe02}) and (\ref{qhe03}), we obtain from (\ref{qhe01})
\begin{align*}
2k\sum_{i=1}^n \sigma_i a(q_{r}^i(\bxi), (\td_h)^{-1}\pt\bxi_i) \le \ve k\sum_{i=1}^n \sigma_i 
\|\pt\bxi_i\|_{-1}^2+K k\sum_{i=1}^n e^{2\alpha t_i}\|\nabla\hbxi_i\|^2.
\end{align*}
Apply (\ref{ieens1}) to have the desired result and this completes the rest of the proof.
\end{proof}
\begin{remark}
The generic error constant $K>0$, involving the estimates of $\bxi_n$ in various norms,
established above, is independent of $n$ and hence the estimates of $\bxi_n$ are uniform in time. In other words, these
estimates are still valid as $t_N\to +\infty$.
\end{remark}
\noindent
We now obtain estimates of $\bta$ below. Hence forward,  $K_n$ means $K(e^{t_n}).$
\begin{lemma}\label{estbta}
Assume (${\bf A1}$), (${\bf A2}$) and a space discretization scheme that satisfies
conditions (${\bf B1}$) and (${\bf B2}$). Further, assume that $\bU^n$ and 
$\bV^n$ satisfy (\ref{fdbej}) and (\ref{fdbejv}), respectively. Then, $\bta_n=\bU^n-\bV^n,~
1\le n\le N$, satisfy the following:
\begin{align}
\|\bta_n\|_r^2+e^{-2\alpha t_n} k\sum_{i=1}^n e^{2\alpha t_i}\|\bta_i\|_{r+1}^2 \le K_n [k(1+\log\frac{1}{k})]^{(1-r)},~~r\in\{0,1\}. \label{estbta1}
\end{align}
\end{lemma}

\begin{proof}
For $n=i,$ we put $\bphi_h=\bta_i$ in (\ref{eebenl}), multiply by $ke^{2\alpha t_i}$
and sum over $1\le i\le n\le N$ to obtain as in (\ref{pree02})
\begin{align}\label{estbta01}
\|\tbta_n\|^2 +2\mu k\sum_{i=1}^n \|\nabla\tbta_i\|^2\le 2k\sum_{i=1}^{n-1} \frac{(e^{2\alpha k}-1)}{k}\|\tbta_i\|^2+2k\sum_{i=1}^n e^{2\alpha t_i} \Lambda_h^i(\bta_i).
\end{align}
We recall from (\ref{dLbe}) that
\begin{align}\label{estbta02}
\Lambda^i_h(\bta_i)= -b(\bu_h^i,\bxi_i,\bta_i)-b(\e_i,\bu_h^i,\bta_i)
-b(\e_i,\bxi_i,\bta_i).
\end{align}
Using (\ref{nonlin1}) and Lemma \ref{pree}, we obtain the following estimates:
\begin{eqnarray}
b(\e_i,\bxi_i,\bta_i)& \le & \|\bxi_i\|^{1/2}\|\nabla\bxi_i\|^{3/2}\|\nabla\bta_i\| +\|\nabla\bxi_i\|\|\bta_i\|\|\nabla\bta_i\| \nonumber \\
& \le & \ve \|\nabla\bta_i\|^2+K\|\bta_i\|^2+Kk^{1/2}(1+\log\frac{1}{k})^{1/2}
\|\nabla\bxi_i\|^2, \label{estbta02a} \\
b(\bu_h^i,\bxi_i,\bta_i) & \le & \ve \|\nabla\bta_i\|^2+K\|\nabla\bxi_i\|^2,
 \label{estbta02b} \\
b(\e_i,\bu_h^i,\bta_i) & \le & \ve \|\nabla\bta_i\|^2+K\big(\|\nabla\bxi_i\|^2
+\|\bta_i\|^2\big). \label{estbta02c}
\end{eqnarray}
Incorporate (\ref{estbta02a})-(\ref{estbta02c}) in (\ref{estbta02}) and then in
(\ref{estbta01})and choose $\ve= \mu/6$ to obtain
\begin{align*}
\|\tbta_n\|^2 +\frac{3}{2}\mu k\sum_{i=1}^n \|\nabla\tbta_i\|^2\le 2k\sum_{i=1}^{n-1} \frac{(e^{2\alpha k}-1)}{k}\|\tbta_i\|^2+K k\sum_{i=1}^n \|\nabla\tbxi_i\|^2 +K k\sum_{i=1}^n \|\tbta_i\|^2.
\end{align*}
Using the boundedness of $\bta_n$, last term on the right hand side can be bounded as 
$$K k\sum_{i=1}^n \|\tbta_i\|^2 = K k \|\tbta_n\|^2 + K k\sum_{i=1}^{n-1}\|\tbta_i\|^2 \le K k e^{2\alpha t_n} + K k\sum_{i=1}^{n-1}\|\tbta_i\|^2.$$
A use of the  Lemma \ref{pree} with the discrete Gronwall's Lemma completes the proof for the case $r=0$. For the case $r=1$, the estimate follows similarly.
\end{proof}



\begin{remark}
Combining Lemmas \ref{pree} and \ref{estbta}, we find suboptimal order of convergence for $\|\e_n\|$:
\begin{eqnarray}
\|\e_n\|^2+e^{-2\alpha t_n} k\sum_{i=1}^n e^{2\alpha t_i}\|\e_i\|_1^2 \le K_nk(1+\log\frac{1}{k}),~1\le n\le N, \label{err1} \\
\|\e_n\|_1^2+e^{-2\alpha t_n} k\sum_{i=1}^n e^{2\alpha t_i}\|\e_i\|_2^2 \le K_n,
~1\le n\le N. \label{err2}
\end{eqnarray}
\end{remark}

\noindent
Below, we shall prove optimal estimate of $\|\e_n\|$ with the help of a series of Lemmas.
\begin{lemma}\label{negbta} 
Under the assumptions of Lemma \ref{estbta}, the following holds:
\begin{align*}
\|\bta_n\|_{-1}^2+e^{-2\alpha t_n} k\sum_{i=1}^n e^{2\alpha t_i} \|\bta_i\|^2
\le K_nk^2(1+\log\frac{1}{k}),~1\le n\le N.
\end{align*}
\end{lemma}

\begin{proof}
Put $\bphi_h= e^{2\alpha t_i} (-\td_h)^{-1}\bta_i$ in (\ref{eebenl}) for $n=i$. Multiply the 
equation by $ke^{2\alpha ik}$ and sum over $1\le i\le n\le N$ to arrive at
\begin{align}\label{negbta01}
\|\tbta_n\|_{-1}^2+2\mu k &\sum_{i=1}^n  \|\tbta_i\|^2 \le \sum_{i=1}^{n-1} 
(e^{2\alpha k}-1)\|\tbta_i\|_{-1}^2+2k \sum_{i=1}^n e^{2\alpha t_i}
\Lambda_h^i (e^{2\alpha t_i} (-\td_h)^{-1}\bta_i ).
\end{align}
From (\ref{dLbe}), we find that
\begin{align}\label{negbta02}
 |2\Lambda_h^i((-\td_h)^{-1}\bta_i)| & \le |2b(\e_i,\bu_h^i,(-\td_h)^{-1}\bta_i) 
+b(\bu_h^i,\e_i,(-\td_h)^{-1}\bta_i) -b(\e_i,\e_i,(-\td_h)^{-1}\bta_i)|.
\end{align}
For the first two terms on the right-hand side of (\ref{negbta02}), we use (\ref{nonlin1}) and similar technique as in (\ref{dth003}) to find that
\begin{equation}\label{negbta03}
|2b(\e_i,\bu_h^i,-\td_h^{-1}\bta_i)|+|2b(\bu_h^i,\e_i,-\td_h^{-1}\bta_i)| \le c\|\e_i\|\|\bu_h^i\|_1\|\bta_i\|_{-1}^{1/2}
\|\bta_i\|^{1/2}
\end{equation}
and for the third terms on the right-hand side of (\ref{negbta02}), we use (\ref{nonlin1}) to observe that
\begin{align}\label{negbta04}
|2b(\e_i,\e_i,-\td_h^{-1}\e_i)| \le c\|\e_i\|\big(\|\e_i\|_1+\|\e_i\|^{1/2}
\|\e_i\|_1^{1/2}\big)\|\bta_i\|_{-1}^{1/2}\|\bta_i\|^{1/2}.
\end{align}
Now, combine (\ref{negbta02})-(\ref{negbta04}) and use the fact that
$ \|\e_i\|_1 \le \|\bu_h^i\|_1+\|\bU^i\|_1 \le K$
to conclude that
\begin{align}\label{negbta07}
 |2\Lambda_h^i((-\td_h)^{-1}\bta_i)| & \le K\|\e_i\|\|\bta_i\|_{-1}^{1/2} 
\|\bta_i\|^{1/2} 
 \le K\|\bxi_i\|\|\bta_i\|_{-1}^{1/2} \|\bta_i\|^{1/2}
 +K\|\bta_i\|_{-1}^{1/2}\|\bta_i\|^{3/2}.
\end{align}
Incorporate (\ref{negbta07}) in (\ref{negbta01}) and use kickback argument to obtain
\begin{align*}
\|\tbta_n\|_{-1}^2+\mu k &\sum_{i=1}^n  \|\tbta_i\|^2 \le Kk\sum_{i=1}^n 
\|\tbta_i\|_{-1}^2+Kk \sum_{i=1}^n \|\tbxi_i\|^2.
\end{align*}
Similar to the proof of Lemma \ref{estbta}, the first term on the right hand side can be written as 
$$K k\sum_{i=1}^n \|\tbta_i\|_{-1}^2 = K k \|\tbta_n\|_{-1}^2 + K k\sum_{i=1}^{n-1}\|\tbta_i\|_{-1}^2 \le K k \|\tbta_n\|^2 + K k\sum_{i=1}^{n-1}\|\tbta_i\|_{-1}^2.$$
An application of Lemmas \ref{neg} and \ref{estbta} with the discrete Gronwall's lemma yields after multiplication by $e^{-2\alpha t_i}$ to the desired estimate. This concludes the rest of the proof.
\end{proof}

\begin{remark}
From Lemmas \ref{neg} and \ref{negbta}, we have the following estimate
\begin{equation}\label{l2ee}
e^{-2\alpha t_n} k\sum_{i=1}^n e^{2\alpha t_i}\|\e_i\|^2 \le K_nk^2(1+\log\frac{1}{k}).
\end{equation}
\end{remark}

We need another estimate of $\bta$ similar to the one in Lemma \ref{ieens} and the proof 
will follow in a similar line. For that purpose, we multiply (\ref{eebenl}) by $k$ and sum over $1$ to $n$ and similar to (\ref{eebelin1i}), we obtain
\begin{align}\label{eebenli}
 (\pt\hbta_n,\bphi_h)+ \mu a(\hbta_n,\bphi_h)+ a(q^n_r(\hbta),\bphi_h) = k\sum_{i=1}^n 
\Lambda_h^i(\bphi_h).
\end{align}

\begin{lemma}\label{ibta}
Under the assumptions of Lemma \ref{estbta}, the following holds:
\begin{align*}
\|\hbta_n\|^2+ e^{-2\alpha t_n} k\sum_{i=1}^n e^{2\alpha t_i} \|\nabla\hbta_i\|^2 \le  K_nk^2(1+\log\frac{1}{k}),~1\le n\le N.
\end{align*}
\end{lemma}

\begin{proof}
Choose $\bphi_h=\hbta_i$ in (\ref{eebenli}) for $n=i$, multiply by $k e^{2\alpha t_i}$ and 
then sum over $1\le i\le n$ to observe as in (\ref{ieens01}) that
\begin{align*}
e^{2\alpha t_n}\|\hbta_n\|^2+\mu k\sum_{i=1}^n e^{2\alpha t_i}\|\nabla\hbta_i\|^2 \le 
k\sum_{i=1}^n e^{2\alpha t_i} k\sum_{j=1}^i |\Lambda_h^i(\hbta_i)|.
\end{align*}
In (\ref{dLbe}), use (\ref{nonlin1}), (\ref{err1}), (\ref{err2}) and (\ref{l2ee}) to obtain 
\begin{align*}
k\sum_{i=1}^n e^{2\alpha t_i}k\sum_{j=1}^i |\Lambda_h^j(\hbta_i)|&= k\sum_{i=1}^n e^{2\alpha t_i} k\sum_{j=1}^i \Big|b(\bu_h^j, \e_j,\hbta_i)+b(\e_j,\bu_h^j,\hbta_i)+b(\e_j,\e_j,\hbta_i)\Big|\nonumber\\
 &\le Kk\sum_{i=1}^n e^{2\alpha t_i} \Big(k\sum_{j=1}^i \big(\|\e_j\| \|\td\bu_h^j\| + \|\e_j\|^{1/2} \|\nabla\e_j\|^{3/2}\big)\Big) \|\nabla\hbta_i\|\nonumber \\
 &\le  K_nk^2(1+\log\frac{1}{k}) +\ve k\sum_{i=1}^n e^{2\alpha t_i}\|\nabla\hbta_i\|^2.
\end{align*}
%
 Now choose $\ve=\mu/2$, we conclude the rest of the proof.
\end{proof}
\noindent We present below a Lemma with optimal estimate for $\bta_n$.
\begin{lemma}\label{eebta} 
Under the assumptions of Lemma \ref{estbta}, the following holds:
\begin{align}\label{eebta1}
t_n\|\bta_n\|^2+e^{-2\alpha t_n} k\sum_{i=1}^n \sigma_i\|\bta_i\|_1^2 \le K_nk^2(1+\log\frac{1}{k}),~1\le n\le N.
\end{align}
\end{lemma}

\begin{proof}
We choose $\bphi_h=\sigma_i\bta_i$ in (\ref{eebenl}) for $n=i$. Multiply the 
resulting equation by $k$ and sum it over $1<i<n$ to find that
\begin{align}\label{eebta01}
\sigma_n\|\bta_n\|^2+2\mu k\sum_{i=1}^n \sigma_i\|\nabla\bta_i\|^2 \le K(\alpha) k\sum_{i=2}^{n-1}\|\tbta_i\|^2-2k\sum_{i=1}^n a(q_{r}^i(\bta),\sigma_i\bta_i) 
+2k\sum_{i=1}^n\Lambda_h^i(\sigma_i\bta_i).
\end{align}
Following the proof of Lemma \ref{quad.hat.est}, we obtain
\begin{align}\label{eebta02}
2k\sum_{i=1}^n \sigma_i a(q_{r}^i(\bta), \bta_i) \le \ve k\sum_{i=1}^n \sigma_i 
\|\nabla\bta_i\|^2+K k\sum_{i=1}^n e^{2\alpha t_i}\|\nabla\hbta_i\|^2.
\end{align}
We recall (\ref{dLbe}) and 
using (\ref{nonlin1}) and similar argument as (\ref{dth003}), to obtain the estimates for nonlinear terms as:
\begin{align}\label{eebta03}
 2k\sum_{i=1}^n\Lambda_h^i(\sigma_i\bta_i) \le \ve k\sum_{i=1}^n \sigma_i\|\nabla\bta_i\|^2
+Kk\sum_{i=1}^n\sigma_i \big(\|\td_h\bu_h^i\|^2+\|\td_h\bU^i\|^2\big)\|\e_i\|^2.
\end{align}
Substitute (\ref{eebta02})-(\ref{eebta03})  in  (\ref{eebta01})  and this completes the rest of the proof.
\end{proof}

\begin{theorem}\label{l2eebe}
 Under the assumptions of Lemma \ref{estbta}, following holds:
\begin{equation}\label{l2eebe1}
 \|\e_n\| \le K_n t_n^{-1/2}k(1+\log \frac{1}{k})^{1/2},~1\le n\le N.
\end{equation}
\end{theorem}

\begin{proof}
Combine the Lemmas \ref{l2eebxi} and \ref{eebta} to complete the rest of the proof.
\end{proof}


\section[Uniform Error]{Uniform Error Estimate}
\se

In this section, we prove the estimate (\ref{l2eebe1}) to be uniform under the uniqueness 
condition $\mu-2N\nu^{-1}\|\f\|_{\infty} >0,$ where $N$ is given as in (\ref{uc}).
We observe that the estimate (\ref{l2eebxi1})
involving $\bxi_n$  is uniform in nature. Hence, we are left to deal with $\bL^2$
estimate of $\bta_n.$

\begin{lemma}\label{btaunif}
Let the assumptions of Lemma \ref{pree} hold. Under the uniqueness condition  $\mu-2N\nu^{-1}\|\f\|_{\infty} >0$, 
the following uniform estimate holds:
\begin{equation*}
 \|\bta_n\| \le K\tau_n^{-1/2}k(1+\log \frac{1}{k})^{1/2},
\end{equation*}
where $\tau_n=\min \{1,t_n\}$.
\end{lemma}

\begin{proof}
Choose $\bphi_h=\bta_i$ in (\ref{eebenl}) for $n=i$. Then from (\ref{dLbe}) and the definition of $N$ (see (\ref{uc})) with the help of (\ref{nonlin1}) and using $\|\e_i\|_2 \le \|\bu_h^i\|_2+\|\bU^i\|_2 \le Kt_i^{-1/2}$, we now bound the nonlinear terms as
%
\begin{align}\label{006}
|\Lambda^i_h(\bta_i)| \le N\|\nabla\bta_i\|^2\|\nabla\bu^i_h\|+K\tau_i^{-3/4}k 
(1+\log\frac{1}{k})^{1/2}\|\nabla\bta_i\|.
\end{align}
We recall from \cite{WHS10} that
$$ \limsup_{t\to+\infty} \|\nabla\bu_h(t)\| \le \nu^{-1}\|\f\|_{\infty}, $$
and therefore, for large enough $i\in \mathbb{N}$, say $i>i_0$ we obtain from (\ref{006})
\begin{align}\label{007}
 |2\Lambda_h^i(\e_i)| \le 2N\nu^{-1}\|\f\|_{\infty}\|\nabla\bta_i\|^2+K\tau_i^{-3/4} 
k (1+\log\frac{1}{k})^{1/2}\|\nabla\bta_i\|.
\end{align}
With $\sigma_i=\tau_i e^{2\alpha t_i},$ we multiply (\ref{eebenl}) with $\bphi_h=\bta_i$ by $k\sigma_i$ and sum 
over $i_0+1$ to $n$ to obtain
\begin{align}\label{008}
k\sum_{i=i_0+1}^n e^{2\alpha t_i} \{\pt\|\bta_i\|^2+2\mu \|\nabla\bta_i\|^2\}
+2k\sum_{i=i_0+1}^n e^{2\alpha t_i} a(q^i_r(\bta),\bta_i)
\le 2k\sum_{i=i_0+1}^n \sigma_i \Lambda^i_h(\bta_i).
\end{align}
Without loss of generality, we can assume that $i_0$ is big enough, so that, by definition
$\tau_i=1$ for $i\ge i_0$ and hence, $\sigma_i=e^{2\alpha t_i}$ for $i_0+1\le i\le n$. With this observation, we rewrite (\ref{008}) as follows:
{\small \begin{align}\label{008a}
k\sum_{i=i_0+1}^n e^{2\alpha t_i} \{\pt\|\bta_i\|^2+2\mu \|\nabla\bta_i\|^2\}
+2k\sum_{i=1}^n e^{2\alpha t_i} a(q^i_r(\bta),\bta_i) 
\le 2k\sum_{i=i_0+1}^n \sigma_i \Lambda^i_h(\bta_i)+2k\sum_{i=1}^{i_0} e^{2\alpha t_i} a(q^i_r(\bta),\bta_i).
\end{align} }
We observe that the last term on the left hand-side of (\ref{008a}) is non-negative and hence, it is dropped. And
\begin{align}\label{009a}
k\sum_{i=i_0+1}^n e^{2\alpha t_i} \{\pt\|\bta_i\|^2\} 
&= k\sum_{i=i_0+1}^n  \Big[\pt \{e^{2\alpha t_i}\|\bta_i\|^2\} - \big(\frac{e^{2\alpha k}-1}{k}\big)e^{2\alpha t_{i-1}}\|\bta_{i-1}\|^2\Big]  \nonumber\\
&= e^{2\alpha t_n}\|\bta_n\|^2 -\sum_{i=i_0+1}^{n-1} e^{2\alpha t_i} (e^{2\alpha k}-1)\|\bta_i\|^2-e^{2\alpha t_{i_0}} \|\bta_{i_0}\|^2.
\end{align}
A use of (\ref{007}) shows
\begin{align}\label{009b}
2k\sum_{i=i_0+1}^n \sigma_i \Lambda^i_h(\bta_i) \le 2k\sum_{i=i_0+1}^n \sigma_i N\nu^{-1} \|\f\|_{\infty}\|\nabla\bta_i\|^2+Kk\sum_{i=i_0+1}^n \tau_i^{1/4} e^{2\alpha t_i}
k(1+\log\frac{1}{k})^{1/2}\|\nabla\bta_i\|.
\end{align}
Substitute (\ref{009a})-(\ref{009b}) in (\ref{008a}) to obtain
\begin{align}\label{009}
e^{2\alpha t_n}\|\bta_n\|^2-\sum_{i=i_0+1}^{n-1} e^{2\alpha t_i} (e^{2\alpha k}-1)\|\bta_i\|^2+\mu k\sum_{i=i_0+1}^n e^{2\alpha t_i}  \|\nabla\bta_i\|^2 
+ k\sum_{i=i_0+1}^n (\mu-2N\nu^{-1}\|\f\|_{\infty}) e^{2\alpha t_i}  \|\nabla\bta_i\|^2 \nonumber \\
\le e^{2\alpha t_{i_0}}\|\bta_{i_0}\|^2+2k\sum_{i=1}^{i_0} e^{2\alpha t_i} q^i_r(\|\nabla\bta_i\|)\|\nabla\bta_i\| 
+Kk^2\sum_{i=i_0+1}^n \tau_i^{1/4} e^{2\alpha t_i} (1+\log\frac{1}{k})^{1/2}\|\nabla\bta_i\|.
\end{align} 
Under the assumption
$ 1+(\frac{\mu\lambda_1}{2})k \ge e^{2\alpha k}, $
which holds for $0 < \alpha < \min \big\{\alpha_0, \delta, \frac{\mu\lambda_1}{2}\big\}$ with $\alpha_0>0$, we find that
\begin{align}\label{010}
\frac{\mu}{2} k\sum_{i=i_0+1}^n e^{2\alpha t_i}  \|\nabla\bta_i\|^2-\sum_{i=i_0+1}^{n-1} 
e^{2\alpha t_i} (e^{2\alpha k}-1)\|\bta_i\|^2 
=k\sum_{i=i_0+1}^{n} \big(\frac{\mu}{2}-\frac{e^{2\alpha k}-1}{k\lambda_1}\big)\sigma_i
\|\nabla\bta_i\|^2 \ge 0.
\end{align}
Due to the uniqueness condition (\ref{uc}), we observe that
\begin{align}\label{011}
 k\sum_{i=i_0+1}^n (\mu-2N\nu^{-1}\|\f\|_{\infty}) e^{2\alpha t_i}  \|\nabla\bta_i\|^2 
\ge 0.
\end{align}
Following the proof techniques of (\ref{pree08})-(\ref{pree09}), we now obtain
\begin{align}\label{012}
2k\sum_{i=1}^{i_0} e^{2\alpha t_i} q^i_r(\|\nabla\bta_i\|) \|\nabla\bta_i\| \le K 
k\sum_{i=1}^{i_0} e^{2\alpha t_i} \|\nabla\bta_i\|^2,
\end{align}
and apply the Young's inequality to arrive at
\begin{align}\label{013}
Kk^2\sum_{i=i_0+1}^n \tau_i^{1/4} e^{2\alpha t_i} (1+\log\frac{1}{k})^{1/2} 
\|\nabla\bta_i\|\le \frac{\mu}{4} k\sum_{i=i_0+1}^n \sigma_i  \|\nabla\bta_i\|^2
+Kk^2 (1+\log\frac{1}{k}) k\sum_{i=i_0+1}^n e^{2\alpha t_i} \tau_i^{-1/2}.
\end{align}
Incorporate (\ref{011})-(\ref{013}) in (\ref{009}), use Lemma \ref{negbta} and
(\ref{l2ee}) to find that
\begin{align*}
e^{2\alpha t_n}\|\bta_n\|^2+k\sum_{i=1}^n \sigma_i  \|\nabla\bta_i\|^2 \le K_{i_0} 
k^2+Kk^2 (1+\log\frac{1}{k}) e^{2\alpha t_n}+Kk\sum_{i=1}^{i_0} e^{2\alpha t_i} \|\nabla\bta_i\|^2.
\end{align*}
Multiply by $e^{-2\alpha t_i}$ and under the assumption, which is proved in the subsequent Lemma \ref{bta1}
\begin{equation}\label{015}
k\sum_{i=1}^{t_0} e^{2\alpha t_i} \|\nabla\bta_i\|^2 \le
K_{i_0}t_{i_0}^{-1}k^2(1+\log \frac{1}{k}),
\end{equation}
we conclude that
$$ \|\bta_n\| \le K_{i_0}t_n^{-1/2}k(1+\log\frac{1}{k})^{1/2}, n>i_0. $$
Combining this with (\ref{eebta1}) for $n\le i_0$, we obtain the desired result,
since $i_0>0$ is fixed. This completes the rest of the proof.
\end{proof}
\noindent We are now left with the proof (\ref{015}).
\begin{lemma}\label{bta1}
Under the assumption of Lemma \ref{btaunif}, the following holds
$$ k\sum_{i=1}^{i_0} e^{2\alpha t_i} \|\nabla\bta_i\|^2 \le 
K_{i_0} t_{i_0}^{-1}k^2(1+\log \frac{1}{k}). $$
\end{lemma}

\begin{proof}
In (\ref{estbta01}), we use
\begin{align*}
\Lambda^i_h(\bta_i) &= -b(\bu^i_h,\e_i,\bta_i)-b(\e_i,\bU^i,\bta_i) 
 \le \frac{\mu}{4}\|\nabla\bta_i\|^2+K\|\e_i\|^2(\|\td_h\bu^i_h\|+\|\td_h\bU^i\|),
\end{align*}
along with Lemma \ref{negbta} and Theorem \ref{l2eebe} to complete the rest of the proof.
\end{proof}
\begin{theorem}\label{unif}
 Under the assumptions of Lemma \ref{btaunif}, following holds:
\begin{equation*}
 \|\e_n\| \le Kt_n^{-1/2}k(1+\log \frac{1}{k})^{1/2},~~1\le n\le N.
\end{equation*}
\end{theorem}

\begin{proof}
Combine the Lemmas \ref{l2eebxi} and \ref{btaunif} to complete the rest of the proof.
\end{proof}

\noindent Now, combine the Theorem \ref{errest}, \ref{l2eebe} and \ref{unif} to conclude the following result:
\begin{theorem}\label{final}
Under the assumptions of Theorem \ref{errest} and Lemma \ref{pree}, following holds:
\begin{equation*}
 \|\bu(t_n)-\bU^n\| \le K_n t_n^{-1/2}\Big(h^2+k(1+\log \frac{1}{k})^{1/2}\Big),~~1\le n\le N.
\end{equation*}
Moreover, under the uniqueness condition (\ref{uc}), the above result is valid uniformly in time.
\end{theorem}

\section{Numerical Experiments}
\begin{figure}[h!] 
\centering
\includegraphics[scale=.38]{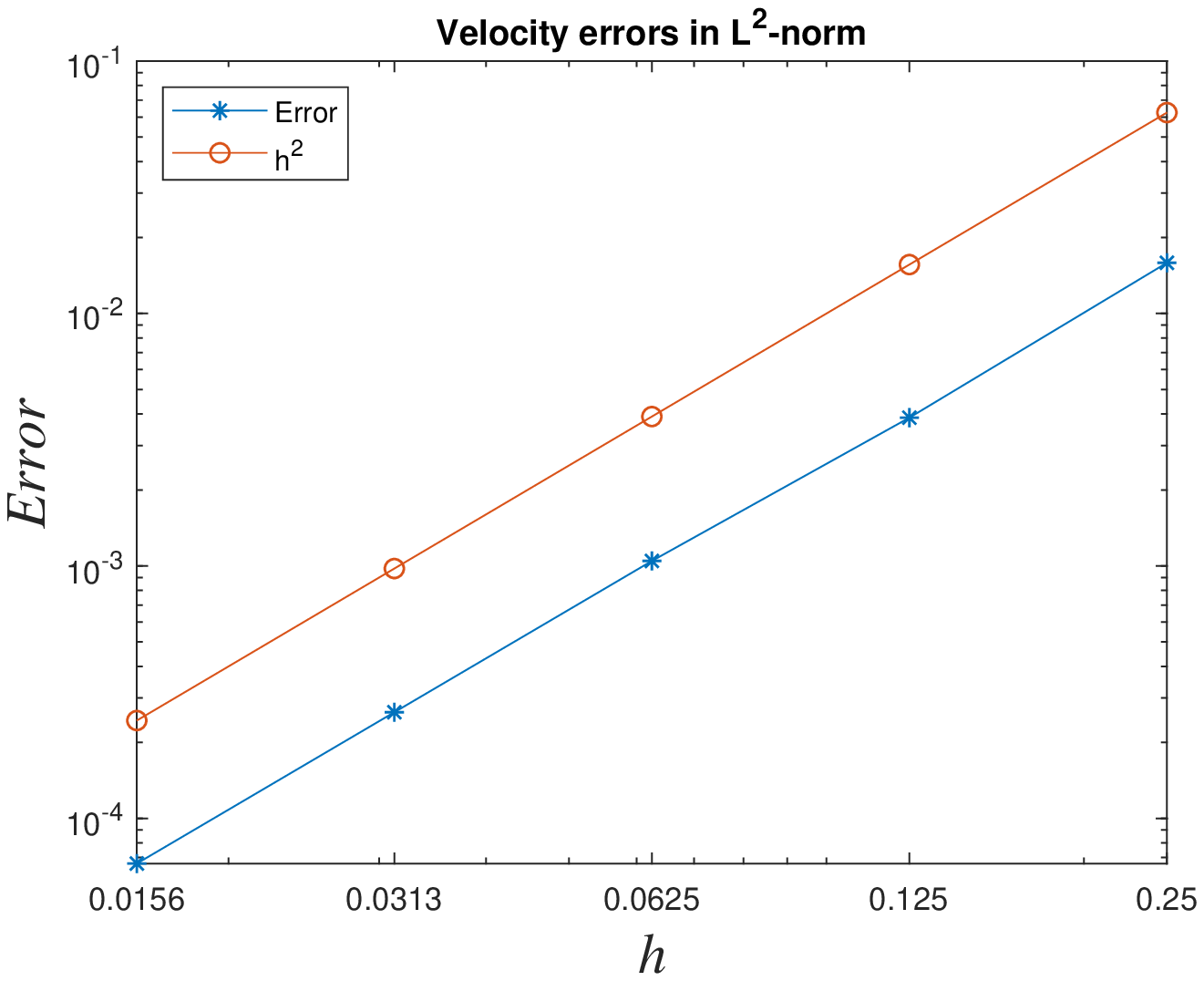}
\includegraphics[scale=.38]{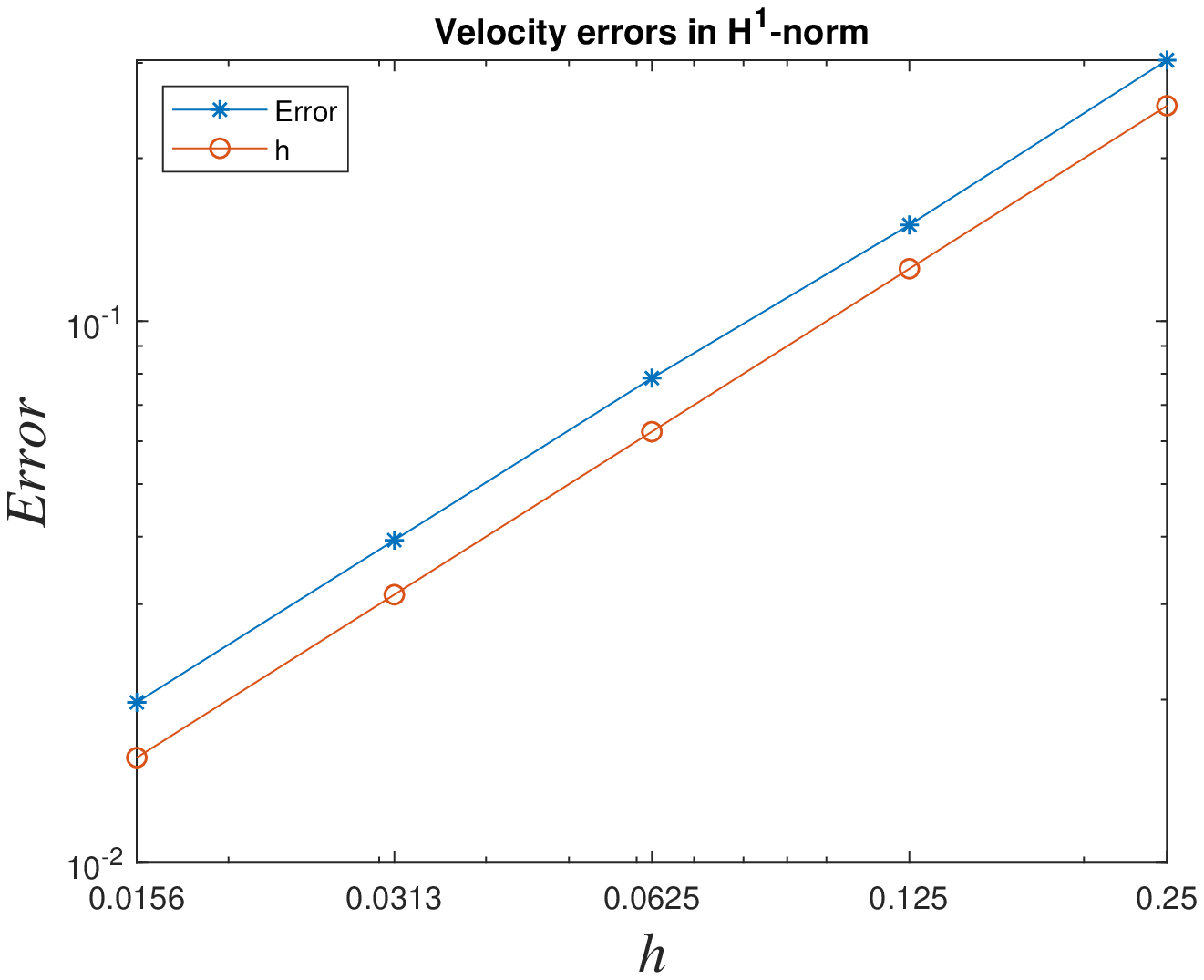}
\includegraphics[scale=.38]{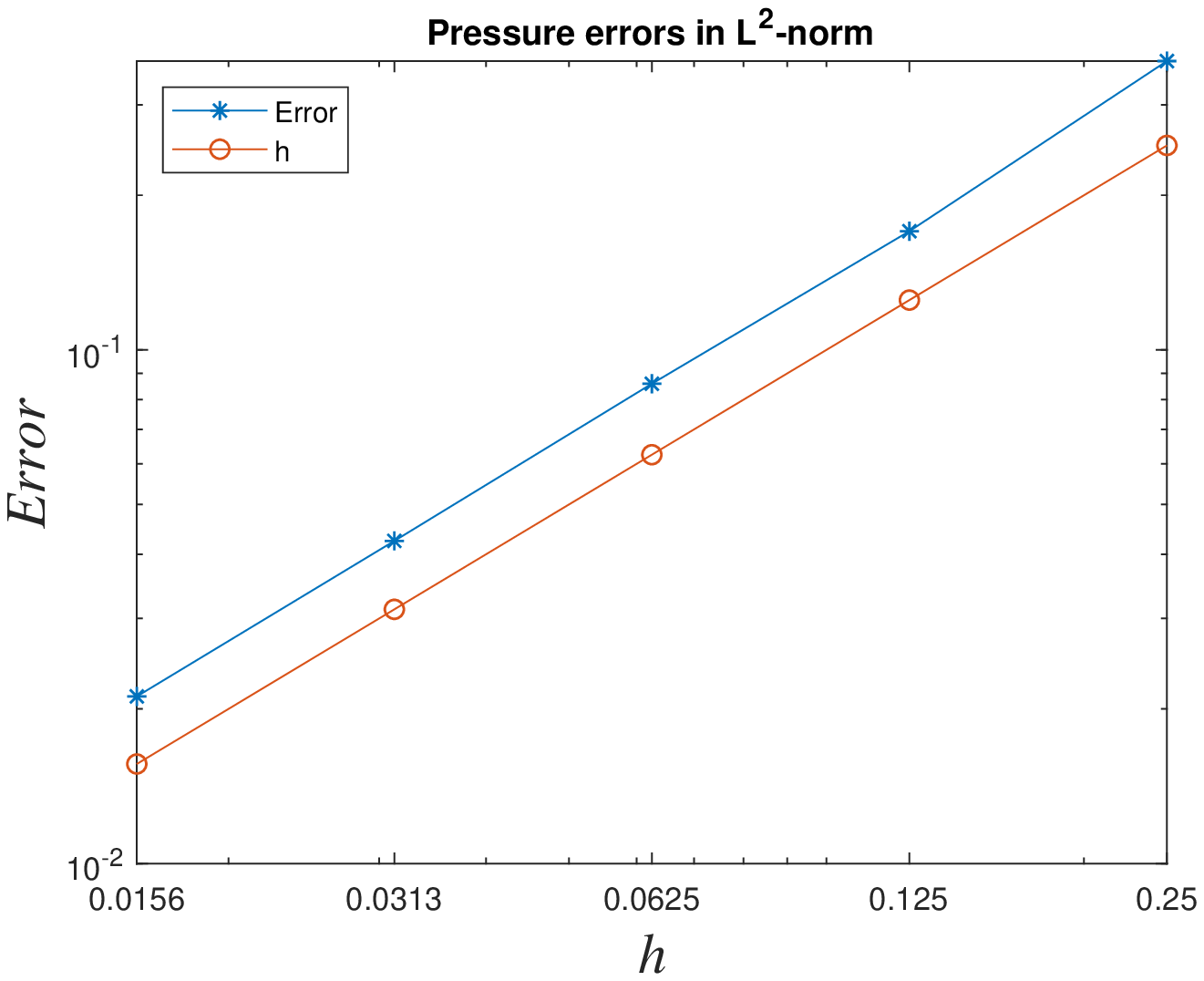}
\caption{Velocity and pressure errors based on P2-P0 element for Example \ref{ex1}.}
\label{fig1}
\end{figure}

\begin{figure}[h!] 
\centering
\includegraphics[scale=.38]{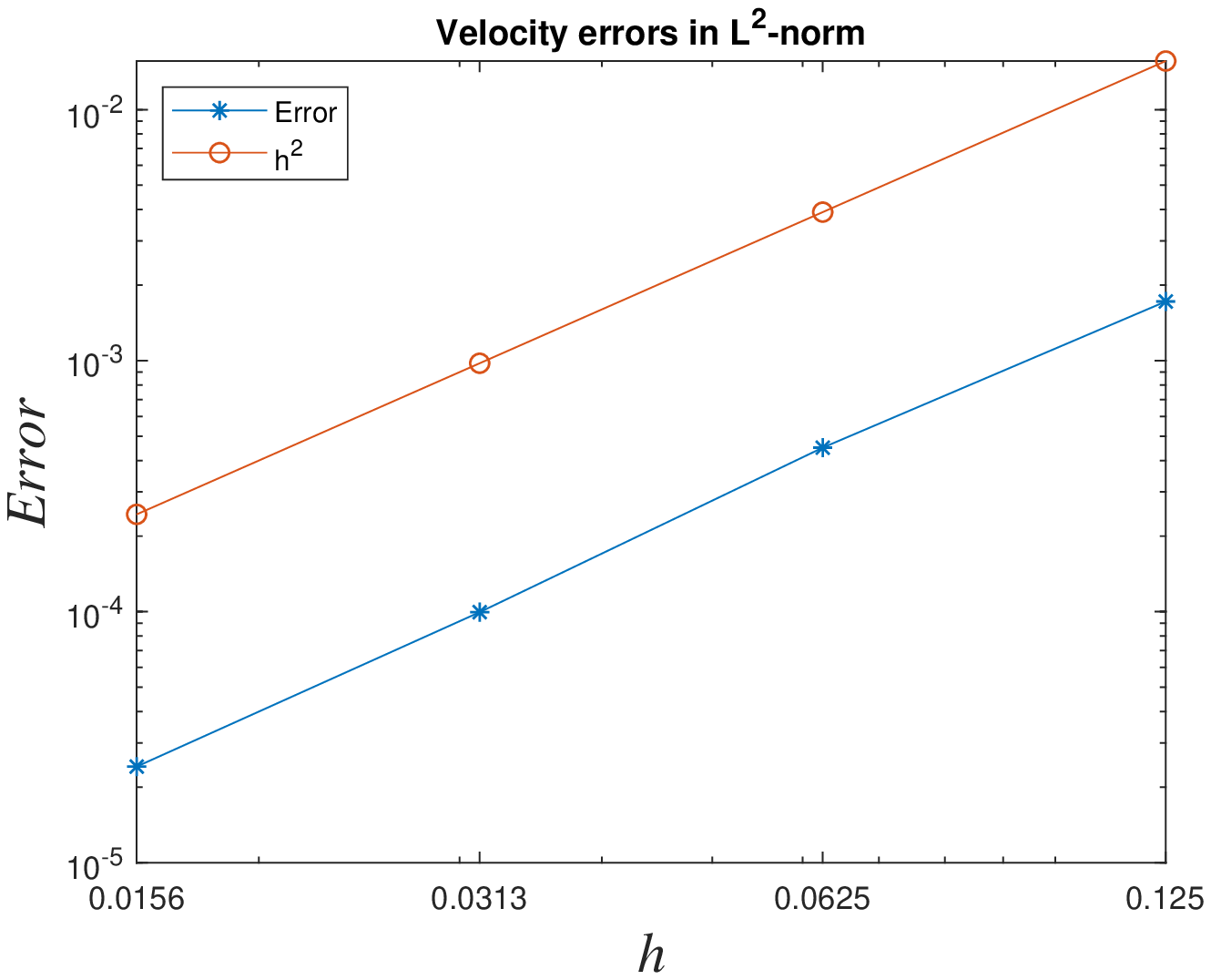}
\includegraphics[scale=.38]{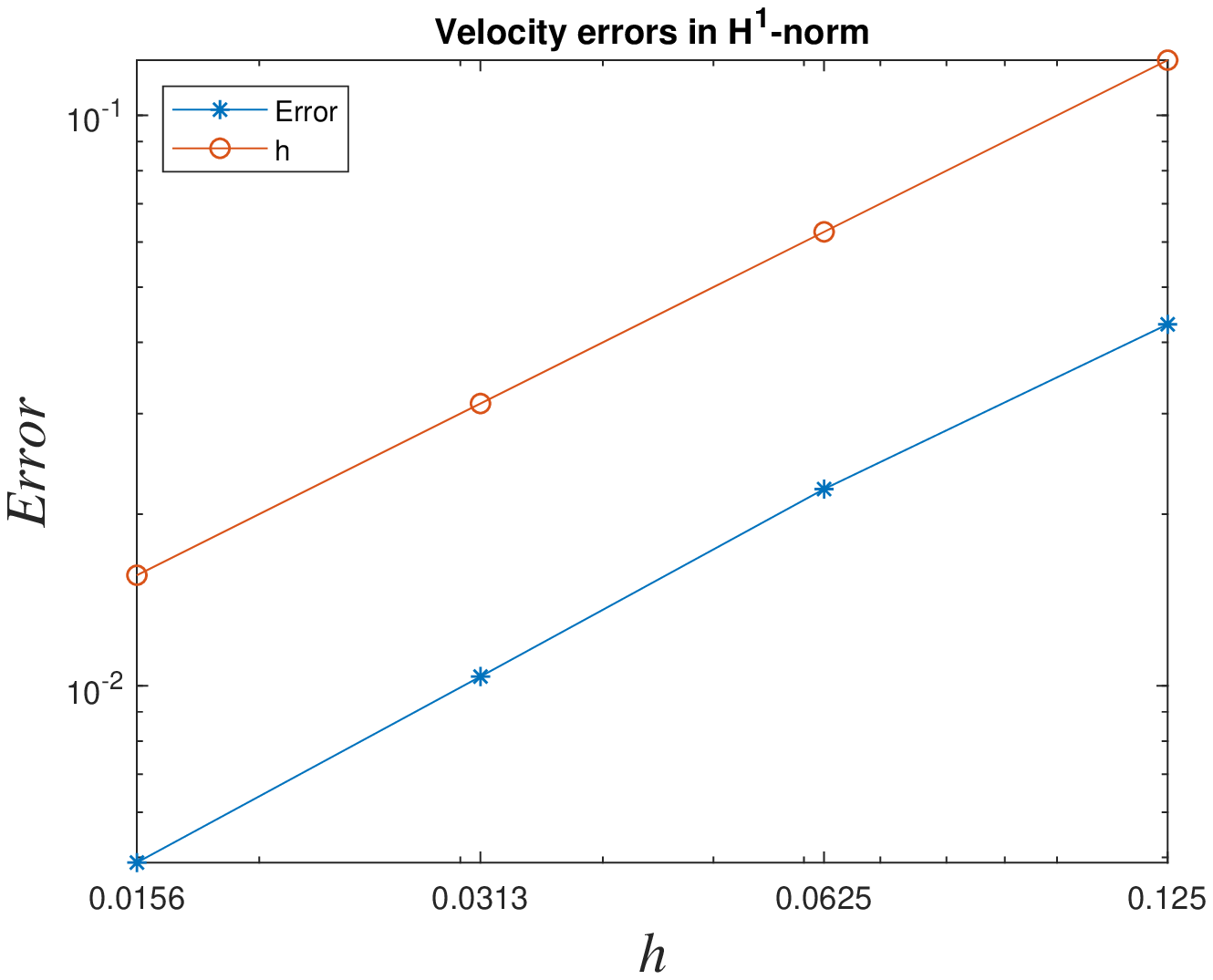}
\includegraphics[scale=.38]{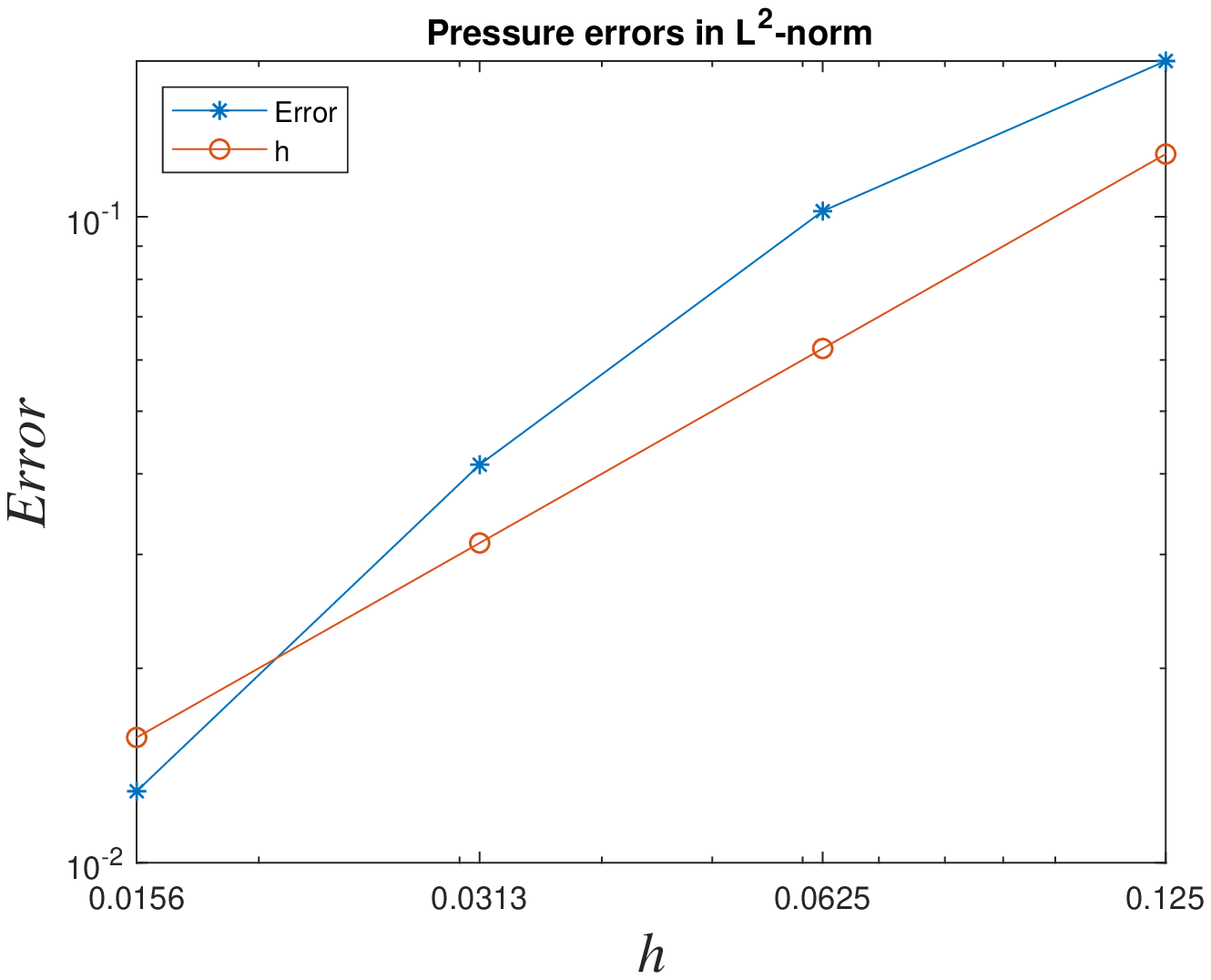}
\caption{Velocity and pressure errors based on MINI element for Example \ref{ex1}.}
\label{fig2}
\end{figure}

In this section, we present some numerical experiments that verify the results of the previous sections,
namely; the order of convergence of the error estimates. For simplicity, we use examples with known solutions. In all cases, computations are done in MATLAB.
%
We consider the Oldroyd model of order one in the domain $\Omega=[0,1]\times[0,1]$ subject to
homogeneous Dirichlet boundary conditions. We approximate the equation using MINI-element and $P2-P0$ element over
a regular triangulation of $\Omega$. The domain is partitioned into triangles with size $h=2^{-i},~i=2,3,\dots,6$. Now, we consider the following example:
\begin{example} \label{ex1}
	For initial data $\bu_0\in\bH^2$, we take the forcing term $f(x,t)$ such that the solution of the problem to be
	\begin{align*}
	u_1(x,t)&= 2e^{t}x^2(x-1)^2y(y-1)(2y-1), \\
	u_2(x,t)&= -2e^{t}x(x-1)(2x-1)y^2(y-1)^2, \\
	p(x,t)&= 2e^t(x-y).
	\end{align*}
\end{example}
\begin{table}[h] 
	\centering %

		\begin{tabular}{|c|c|c|c|c|c|c|c|}
		\hline
		h   & $\|u(t_n)-U^n\|_{L^2}$ & Rate & $\|u(t_n)-U^n\|_{H^1}$ & Rate & $\|p(t_n)-P^n\|_{L^2}$ & Rate   \\
		\hline 
		1/8   &  0.00386700  &         &  0.15057567  &         &  0.17021691  &        \\
		1/16  &  0.00104657  & 1.8855  &  0.07849371  & 0.9398  &  0.08591565  & 0.9864  \\
		1/32  &  0.00026335  & 1.9906  &  0.03939885  & 0.9944  &  0.04246851  & 1.0165  \\
		1/64  &  0.00006623  & 1.9913  &  0.01976541  & 0.9952  &  0.02115282  & 1.0055  \\
		\hline
	\end{tabular}\label{t1}
	\caption{Errors and convergence rates for backward Euler method for Example \ref{ex1} for P2-P0 element}
\end{table}
\begin{table}[h] 
	\centering %
		\begin{tabular}{|c|c|c|c|c|c|c|c|}
		\hline
		h   & $\|u(t_n)-U^n\|_{L^2}$ & Rate & $\|u(t_n)-U^n\|_{H^1}$ & Rate & $\|p(t_n)-P^n\|_{L^2}$ & Rate   \\
		\hline 
		1/8   &  0.00172068  &         &  0.04302980  &         &  0.17416894  &         \\
		1/16  &  0.00045020  & 1.9344  &  0.02212674  & 0.9595  &  0.10199069  & 0.7720  \\
		1/32  &  0.00009954  & 2.1771  &  0.01037882  & 1.0921  &  0.04131507  & 1.3037  \\
		1/64  &  0.00002414  & 2.0436  &  0.00489803  & 1.0834  &  0.01289942  & 1.6794  \\
		\hline
	\end{tabular}\label{t2}
		\caption{Errors and convergence rates for backward Euler method for Example \ref{ex1} for MINI-element}
\end{table}
\noindent
In Table 1 and 2, we present the numerical errors and convergence rates obtained on successive meshes using $P2-P0$ element and MINI-element, respectively, for backward Euler scheme applied to the system (\ref{om})-(\ref{ibc}) with $\mu=1, \gamma=0.1, \delta=0.1$ and final time $T=1$. The theoretical analysis shows that the rate of convergence are $\mathcal{O}(h^2)$ in $\bL^2$-norm and  $\mathcal{O}(h)$ in $\bH^1$-norm for the velocity and $\mathcal{O}(h)$ in $\bL^2$-norm for the pressure with the choice of $k=\mathcal{O}(h^2)$. The error graphs are presented in Fig \ref{fig1} and Fig \ref{fig2}. These results support the optimal convergence rates obtained in Theorem \ref{final}. 


In order to verify the rate of convergence in both spatial and temporal directions and the uniform convergence in time for nonsmooth data, we consider the following example \cite{HHF15, ZQ18}. 

\begin{figure}[h!] 
\centering
\includegraphics[scale=.38]{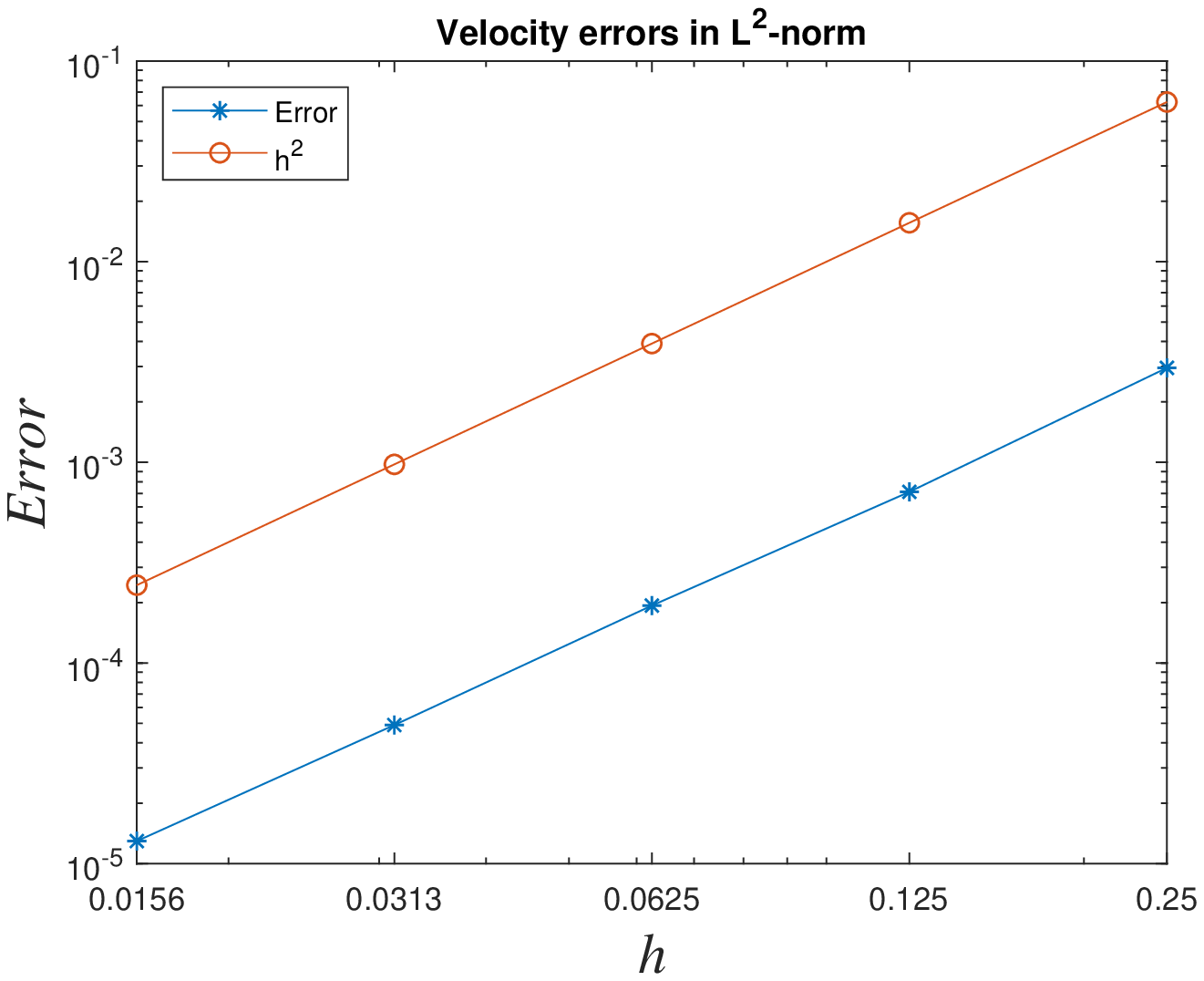}
\includegraphics[scale=.38]{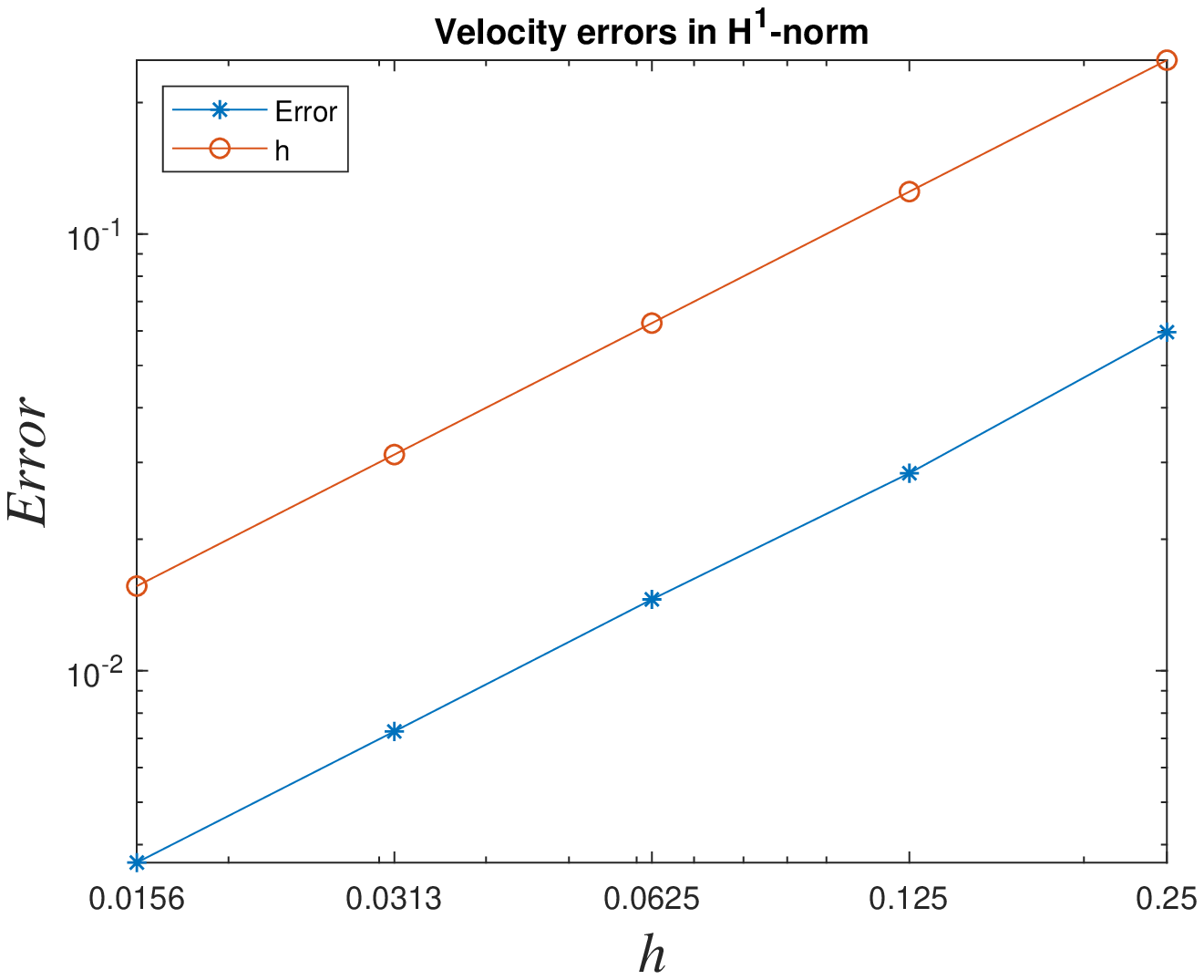}
\includegraphics[scale=.38]{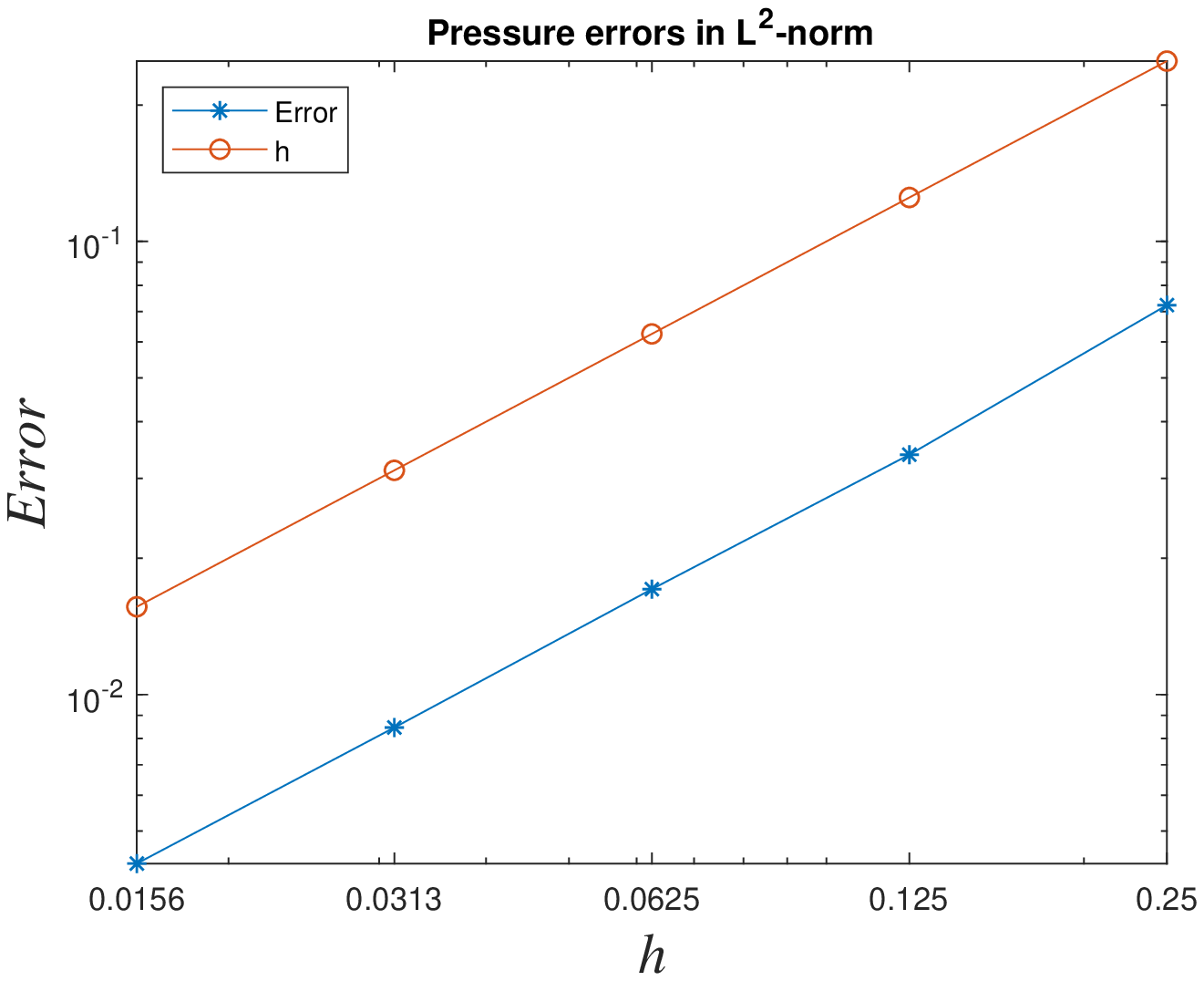}
\caption{Velocity and pressure errors based on P2-P0 element for Example \ref{ex2}.}
\label{fig3}
\end{figure}

\begin{figure}[h!] 
\centering
\includegraphics[scale=.38]{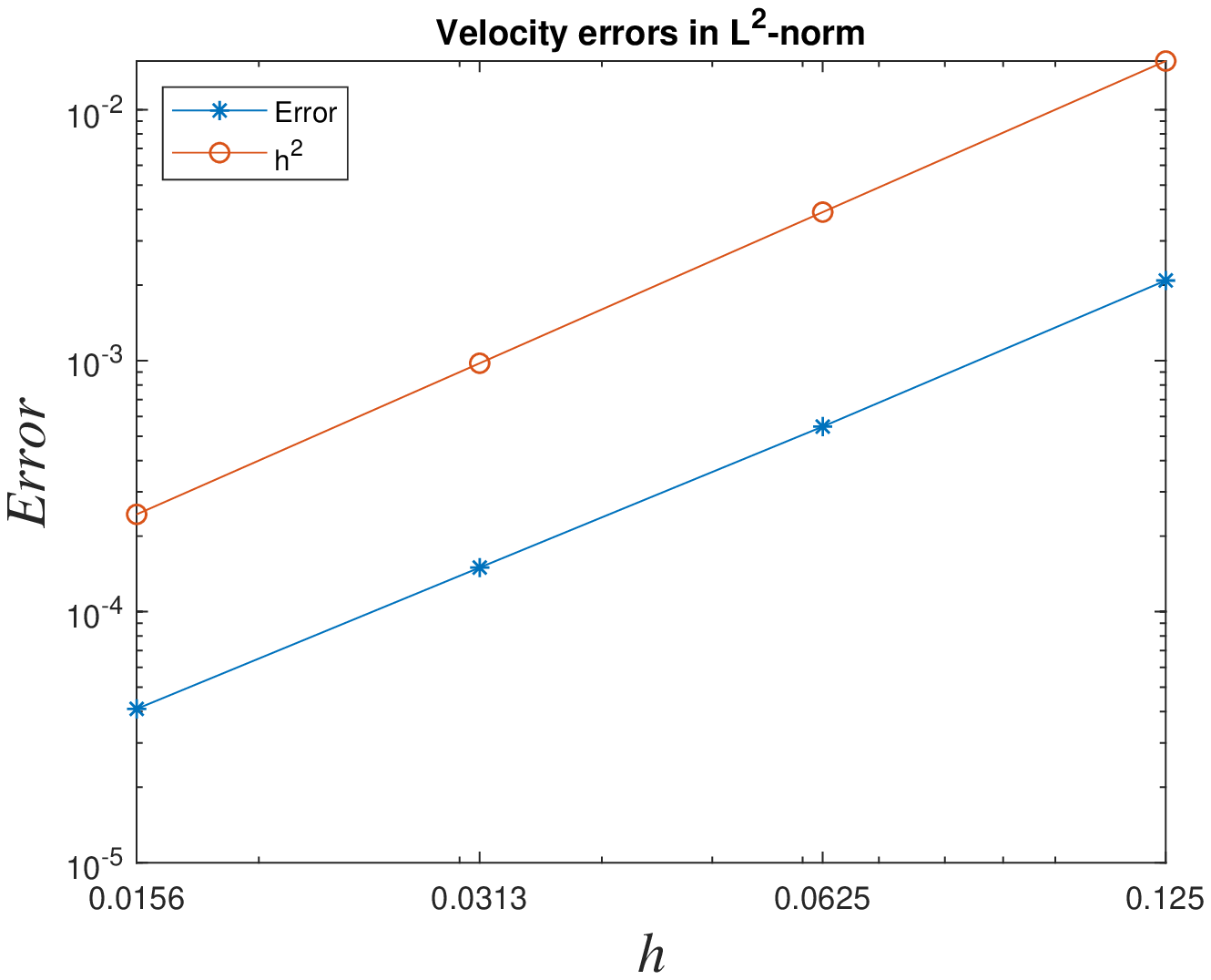}
\includegraphics[scale=.38]{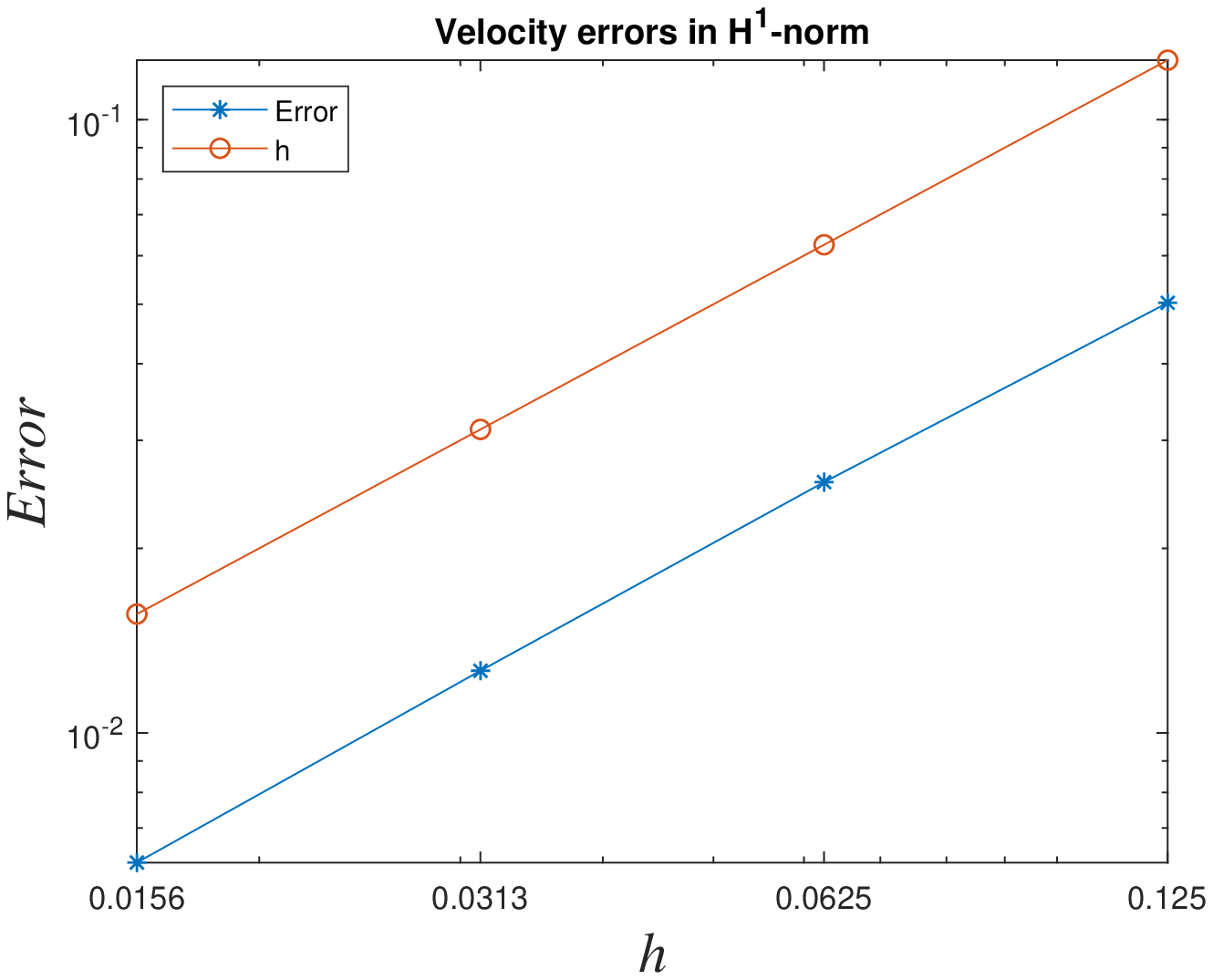}
\includegraphics[scale=.38]{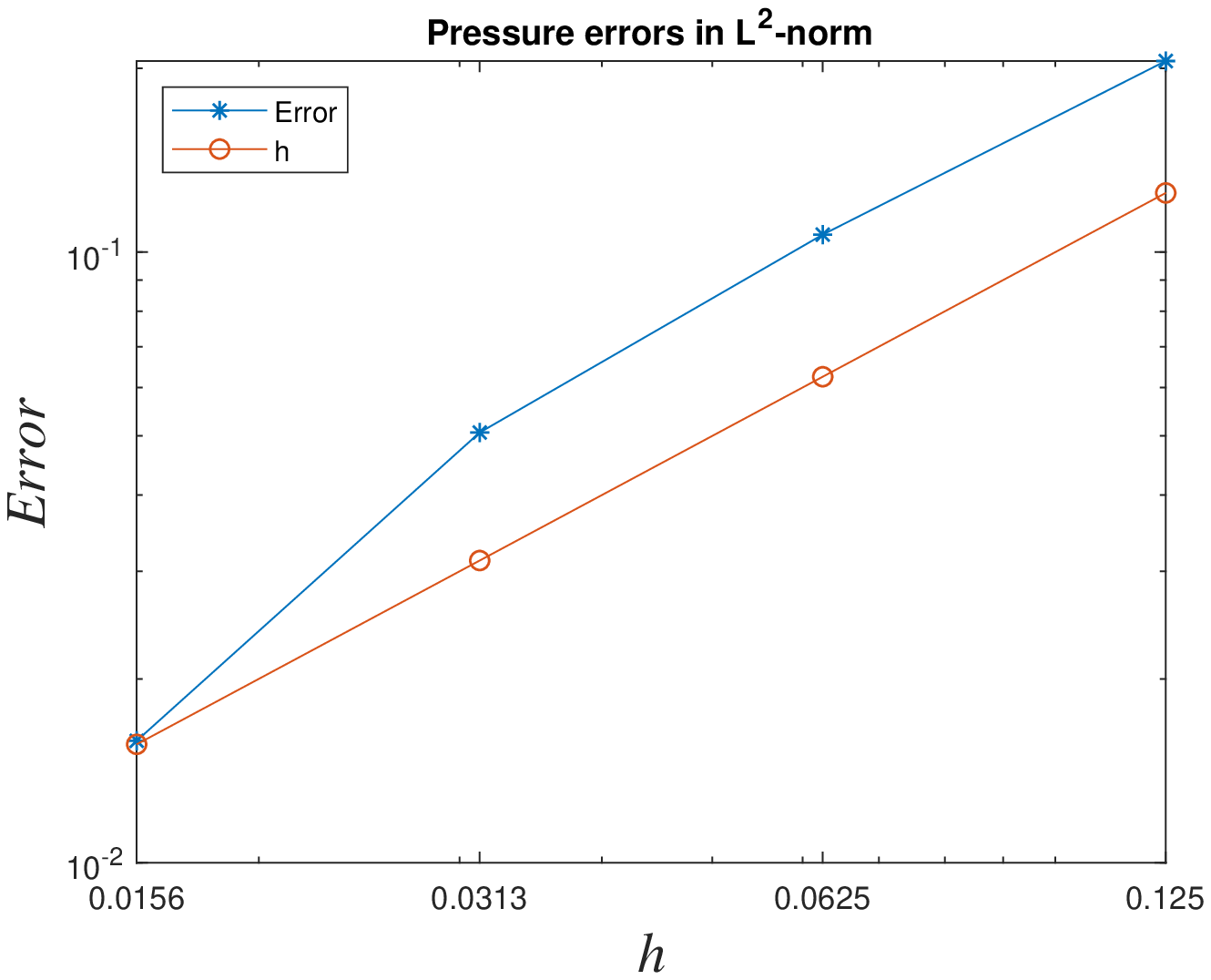}
\caption{Velocity and pressure errors based on MINI element for Example \ref{ex2}.}
\label{fig4}
\end{figure}

\begin{figure}[h!] 
\centering
\includegraphics[scale=.47]{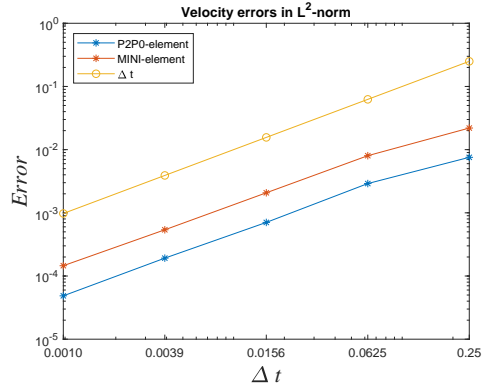}
\caption{Velocity errors in $L^2$- norm with respect to time for Example \ref{ex2}.}
\label{fig5}
\end{figure}

\begin{example}\label{ex2}
For initial data $\bu_0\in\bH_0^1$, we take the forcing term $f(x,t)$ such that the solution of the problem to be
	\begin{align*}
	u_1(x,t)&= 5 x^{5/2}(x-1)^2y^{3/2}(y-1)(9y-5)\cos{t}, \\
	u_2(x,t)&= -5 x^{3/2}(x-1)(9x-5)y^{5/2}(y-1)^2 \cos{t}, \\
	p(x,t)&= 2 (x-y)\cos{t}.
	\end{align*}
\end{example}
\begin{table}[h] 
\centering
\begin{tabular}{|c|c|c|c|c|c|c|c|}
   \hline
   h   & $\|u(t_n)-U^n\|_{L^2}$ & Rate & $\|u(t_n)-U^n\|_{H^1}$ & Rate & $\|p(t_n)-P^n\|_{L^2}$ & Rate  \\
   \hline 
  1/4   &  0.00295597  &         &  0.05958679  &         &  0.07233700  &         \\
  1/8   &  0.00071240  & 2.0529  &  0.02832958  & 1.0727  &  0.03383893  & 1.0960  \\
  1/16  &  0.00019314  & 1.8830  &  0.01456592  & 0.9597  &  0.01708781  & 0.9857  \\
  1/32  &  0.00004903  & 1.9780  &  0.00726227  & 1.0041  &  0.00845973  & 1.0143   \\
  1/64  &  0.00001294  & 1.9217  &  0.00363780  & 0.9973  &  0.00423842  & 0.9971   \\
   \hline
\end{tabular}\label{t3}
\caption{Errors and convergence rates for backward Euler method for Example \ref{ex2} for P2-P0 element}
\end{table}				
\begin{table}[h!] 
	\centering %
		\begin{tabular}{|c|c|c|c|c|c|c|c|}
		\hline
		h   & $\|u(t_n)-U^n\|_{L^2}$ & Rate & $\|u(t_n)-U^n\|_{H^1}$ & Rate & $\|p(t_n)-P^n\|_{L^2}$ & Rate   \\
		\hline 
		1/8   &  0.00208654  &         &  0.05026012  &         &  0.20559044  &         \\
		1/16  &  0.00054627  & 1.9334  &  0.02563557  & 0.9713  &  0.10681639  & 0.9446  \\
		1/32  &  0.00014985  & 1.8660  &  0.01262894  & 1.0214  &  0.05064130  & 1.0767  \\
		1/64  &  0.00004098  & 1.8705  &  0.00614613  & 1.0390  &  0.01583420  & 1.6772  \\
		\hline
	\end{tabular}\label{t4}
		\caption{Errors and convergence rates for backward Euler method for Example \ref{ex2} for MINI-element}
\end{table}
\begin{table}[h!] 
	\centering %
		\begin{tabular}{|c|c|c|c|c|}
		\hline
				& \multicolumn{2}{|c|}{ \textbf{P2-P0 element}} & \multicolumn{2}{|c|}{ \textbf{MINI element}}\\
		\hline
		$k$   & $\|u(t_n)-U^n\|_{L^2}$ & Rate & $\|u(t_n)-U^n\|_{L^2}$ & Rate   \\
		\hline 
		1/4     & 0.00755768	& 			&  0.02197941  & 		 \\
		1/16    & 0.00290626	& 0.6894	&  0.00802012  & 0.7272   \\
		1/64    & 0.00070503	& 1.0217	&  0.00207538  & 0.9751    \\
		1/256   & 0.00019182	& 0.9389	&  0.00053884  & 0.9727    \\
		1/1024  & 0.00004856	& 0.9909	&  0.00014555  & 0.9442   \\
		\hline 
	\end{tabular} \label{t5}
		\caption{$L^2$-Errors and convergence rates in temporal direction for Example \ref{ex2} for P2-P0 and MINI-elements}
\end{table}			
In Table 3 and 4, we have shown the errors and the convergence rates for the backward Euler method using P2-P0 and MINI elements, respectively,
with $\mu=1, \gamma=0.1, \delta=1$ and final time $T=1$. The numerical results confirm the optimal convergence rates of the velocity error in $\bL^2$-norm as in Theorem \ref{final}. The error graphs are given in Fig \ref{fig3} and Fig \ref{fig4}.
In Table 5, we present the errors and the convergence rate in temporal direction for P2-P0 and MINI-elements, respectively. Here, we take $k= 2^{-2i},~i=1,2,\dots,5$, $\mu=1, \delta=0.1, \gamma=0.1, h=\mathcal{O}(\sqrt{k})$ and $T=1$. The error graph is given in Fig \ref{fig5}. We observe that the rate of convergence confirms the theoretical findings.
\begin{table}[h!] 
	\centering %
		\begin{tabular}{|c|c|c|c|c|c|}
		\hline
			& 			& \multicolumn{2}{|c|}{ \textbf{P2-P0 element}} & \multicolumn{2}{|c|}{ \textbf{MINI element}}\\
		\hline
		Final time & $h$   & $\|u(t_n)-U^n\|_{L^2}$ & Rate & $\|u(t_n)-U^n\|_{L^2}$ & Rate   \\
		\hline
		T=10&	1/4   & 0.00475212	& 			&  0.03323502  & 		 \\
			&	1/8   & 0.00116079	& 2.0334	&  0.00317788  & 3.3866   \\
			&	1/16  & 0.00031742	& 1.8706	&  0.00081112  & 1.9701    \\
			&	1/32  & 0.00008496	& 1.9015	&  0.00020500  & 1.9843    \\
			&	1/64  & 0.00002272	& 1.9030	&  0.00005176  & 1.9857   \\
		\hline 
		T=20&	1/4   & 0.00195842	& 			&  0.01616627  & 		 \\
			&	1/8   & 0.00047565	& 2.0417	&  0.00156824  & 3.3657   \\
			&	1/16  & 0.00013005	& 1.8708	&  0.00040039  & 1.9697    \\
			&	1/32  & 0.00003459	& 1.9105	&  0.00010100  & 1.9870    \\
			&	1/64  & 0.00000919	& 1.9120	&  0.00002545  & 1.9885   \\
		\hline
		T=30&	1/4   & 0.00149863	& 			&  0.00610584  & 		 \\
			&	1/8   & 0.00036956	& 2.0197	&  0.00055113  & 3.4697   \\
			&	1/16  & 0.00010101	& 1.8713	&  0.00014395  & 1.9368    \\
			&	1/32  & 0.00002742	& 1.8810	&  0.00004044  & 1.8318    \\
			&	1/64  & 0.00000744	& 1.8824	&  0.00001135  & 1.8333   \\
		\hline 
		T=40&	1/4   & 0.00446125	& 			&  0.02641260  & 		 \\
			&	1/8   & 0.00109409	& 2.0277	&  0.00249141  & 3.4062   \\
			&	1/16  & 0.00029908	& 1.8711	&  0.00064118  & 1.9582    \\
			&	1/32  & 0.00008028	& 1.8973	&  0.00016808  & 1.9316    \\
			&	1/64  & 0.00002153	& 1.8988	&  0.00004402  & 1.9330   \\
		\hline
		T=50&	1/4   & 0.00599158	& 			&  0.03821832  & 		 \\
			&	1/8   & 0.00146697	& 2.0301	&  0.00363026  & 3.3961   \\
			&	1/16  & 0.00040103	& 1.8710	&  0.00093218  & 1.9614    \\
			&	1/32  & 0.00010770	& 1.8967	&  0.00024209  & 1.9451    \\
			&	1/64  & 0.00002889	& 1.8982	&  0.00006281  & 1.9465  \\
		\hline
	\end{tabular} \label{t6}
		\caption{$L^2$-Errors and convergence rates for Example \ref{ex2} for P2-P0 and MINI-elements}
\end{table}
\begin{figure}[h!] 
\centering
\includegraphics[scale=.5]{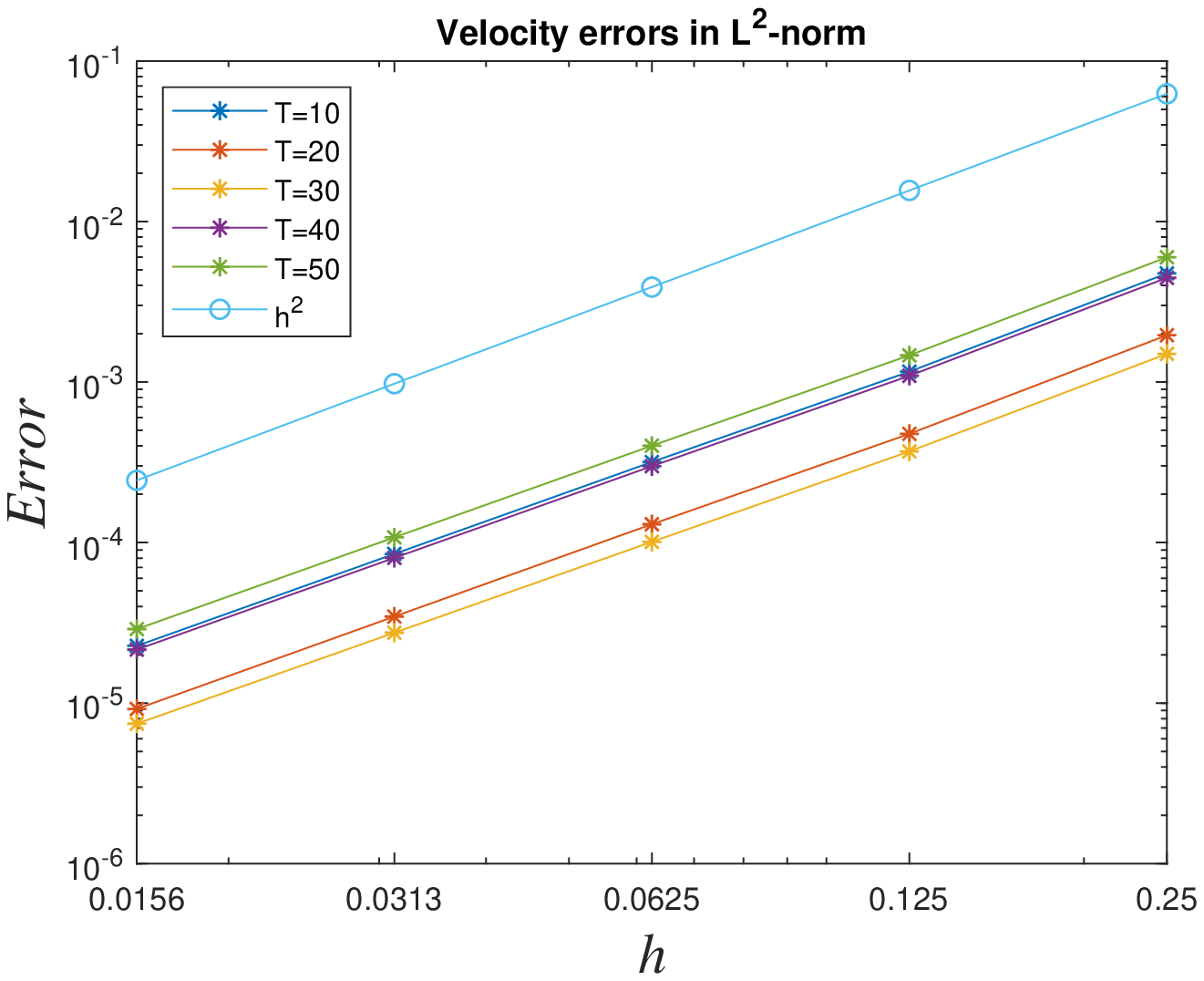}
\includegraphics[scale=.5]{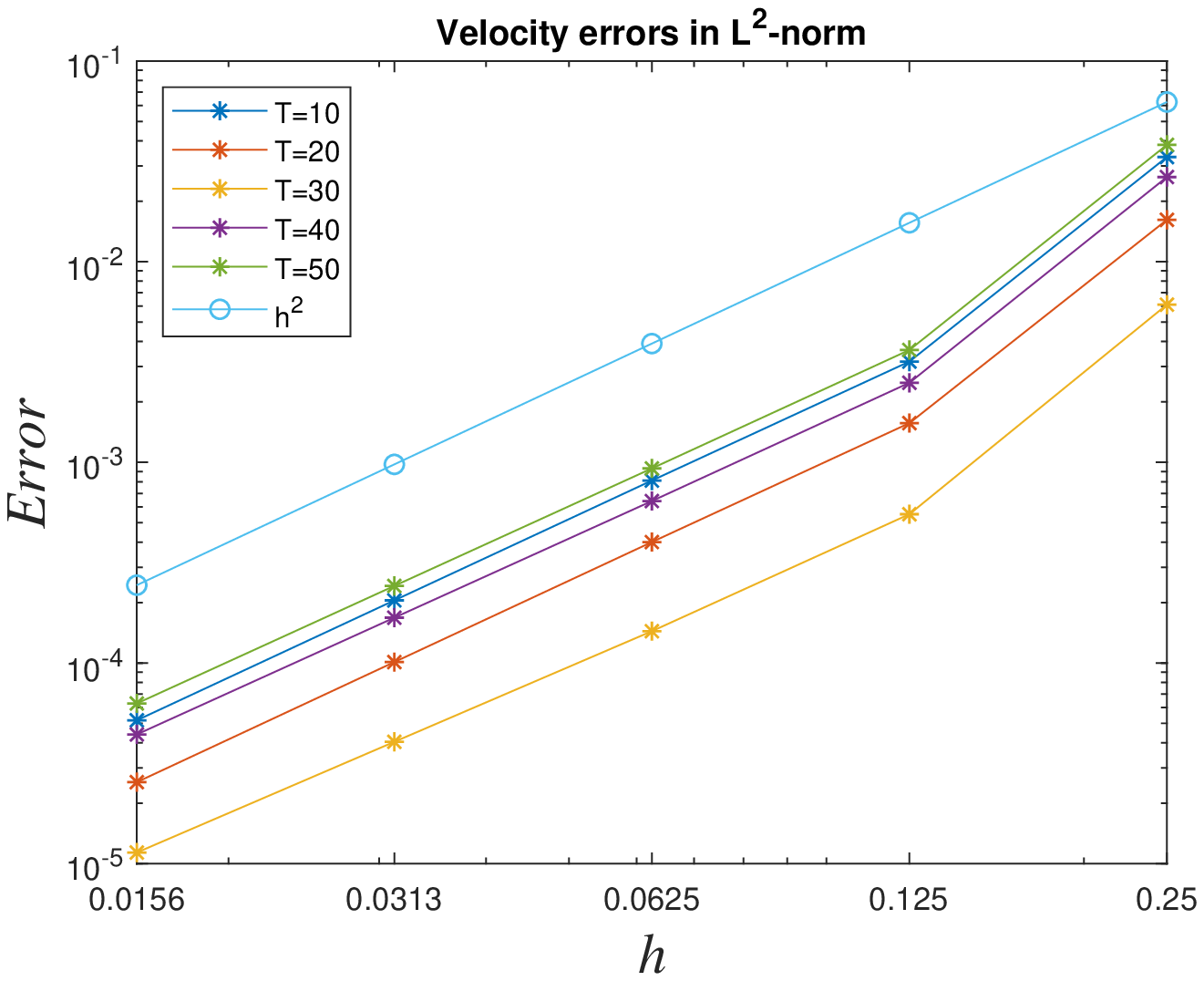}
\caption{Uniform in time errors for P2-P0 element (left) and MINI element (right) for Example \ref{ex2}.}
\label{fig6}
\end{figure}

For the example \ref{ex2}, the numerical results are shown for final time $T=10, 20, 30, 40$ and $50$ with $\mu=1, \gamma=0.1, \delta=1$, $k=0.1$ and $h=2^{-i},~i=2,3,\dots,6$. 
We represent the errors and the convergence rates for the velocity in $\bL^2$-norm for P2-P0 and MINI-elements in Table 6 and Fig \ref{fig6}. The numerical experiments show that for a large time the convergence rates remain same.
%
%
%
%
%
%
%
%
%
%
%
%
%

\begin{table}[h!]\label{t8}
\centering
\begin{tabular}{|c|c|c|c|c|c|}
   \hline
  	& \backslashbox{h}{k}   & 0.1 & 0.5 & 1 & 1.3   \\
   \hline 
  		& 1/10   &  0.04066058  & 0.04013758  &  0.03994627  & 0.03431507    \\
P2-P0   & 1/20   &  0.04060207  & 0.04007989  &  0.03988995  & 0.03426480    \\
element & 1/30   &  0.04059567  & 0.04007359  &  0.03988379  & 0.03425928   \\
  		& 1/40	 &	0.04059327	& 0.04007121  &  0.03988147  & 0.03425720   \\
  \hline 
  		& 1/10   &  0.04164312  & 0.03825401  &  0.03808147  & 0.03272622    \\
MINI    & 1/20   &  0.04310498  & 0.03967860  &  0.03949458  & 0.03393053    \\
element & 1/30   &  0.04334520  & 0.03991041  &  0.03972447  & 0.03412640   \\
 		& 1/40	 &	0.04343595	& 0.03999863  &  0.03981201  & 0.03420103   \\
  \hline
\end{tabular}
\caption{The norm $\sup_{0\le t_n\le 5}\|U^n\|_{L^2}$  with nonsmooth data for P2-P0 and MINI elements}
\end{table}
\begin{table}[h!]\label{t9}
\centering
\begin{tabular}{|c|c|c|c|c|c|}
   \hline
  & \backslashbox{h}{k}   & 0.1 & 0.5 & 1 & 1.3   \\
   \hline 
  		&	1/10   &  0.31297840  & 0.30901528  &  0.30752906  & 0.26430325    \\
P2-P0   &	1/20   &  0.31034508  & 0.30641145  &  0.30500030  & 0.26201939    \\
element &   1/30   &  0.30987407  & 0.30594553  &  0.30454751  & 0.26161070   \\
  		&	1/40   &	0.30970478	& 0.30577810  &  0.30438480  & 0.26146391   \\
  \hline
  		&	1/10   &  0.32158956  & 0.29632714  &  0.29503887  & 0.25354515    \\
MINI  	&	1/20   &  0.32860122  & 0.30352433  &  0.30217202  & 0.25959799    \\
element &  1/30   &  0.32979560  & 0.30474032  &  0.30337746  & 0.26062074   \\
  		&	1/40	 &	0.33024171	& 0.30519856  &  0.30383212  & 0.26100694   \\
  \hline
\end{tabular}
\caption{The norm $\sup_{0\le t_n\le 5}\|U^n\|_{H^1}$  with nonsmooth data for P2-P0 and MINI elements}
\end{table}
\begin{table}[h!]\label{t10}
\centering
\begin{tabular}{|c|c|c|c|c|c|}
   \hline
  & \backslashbox{h}{k}   & 0.1 & 0.5 & 1 & 1.3   \\
   \hline 
  		&	1/10   &  0.80976358  & 0.80204569  &  0.80174532  & 0.69388452    \\
P2-P0   &	1/20   &  0.81411956  & 0.80642016  &  0.80616249  & 0.69778546    \\
element &  1/30   &  0.81499218  & 0.80729509  &  0.80704511  & 0.69856265   \\
  		&	1/40	 &	0.81529345	& 0.80759726  &  0.80735004  & 0.69883114   \\
  \hline
  		&	1/10   &  0.86979545  & 0.85337708  &  0.85455721  & 0.73827331    \\
MINI     & 1/20   &  0.82750272  & 0.81637519  &  0.81639044  & 0.70638134    \\
element &  1/30   &  0.82289194  & 0.81236294  &  0.81225804  & 0.70293783   \\
  		&	1/40	 &	0.81951684	& 0.80935013  &  0.80914669  & 0.70034738   \\
  \hline
\end{tabular}
\caption{The norm $\sup_{0\le t_n\le 5}\|P^n\|_{L^2}$ of  with nonsmooth data for P2-P0 and MINI elements}
\end{table}

In tables 7 to 9, we have shown the  maximal $\bL^2$- and $\bH^1$-norm of the velocity and maximal $L^2$-norm of the pressure among several time  steps $k=0.1, 0.5, 1, 1.3$ again for the Example \ref{ex2}. The results indicate that the scheme can run well for the values of the time steps going from $k = 0.1$ to $k = 1.3$, but there is a deterioration of the convergence rate for $k=1$ and $k=1.3$.


\section{Conclusion}
In this article, optimal error estimates are derived for the backward Euler method employed to the Oldroyd model with non-smooth initial data, that is, $\bu_0\in H_0^1(\Omega)$. For the complete discrete scheme, uniform a priori bounds are shown for the discrete solution. Both optimal and uniform error estimate for the velocity are proved. Uniform estimates are derived under the uniqueness condition. The analysis has been done for the non-smooth initial data and the proofs are more involved in comparison to the smooth case. Our numerical results confirms our theoretical results.

\end{document}